%% file: main.tex
\documentclass[review, 3p, times]{elsarticle}

\input{mathhead}
\begin{document}

\begin{frontmatter}

\title{Enhancing Sparsity of Hermite Polynomial Expansions by Iterative Rotations}

\author[pnnl1]{Xiu Yang}
\author[pnnl1]{Huan Lei}
\author[pnnl2]{Nathan A.~Baker}
\author[purdue1,purdue2]{Guang Lin\corref{gl}}

\address[pnnl1]{Advanced Computing, Mathematics, and Data Division, Pacific
                Northwest National Laboratory, Richland, WA 99352, USA}
\address[pnnl2]{Computational and Statistical Analytics Division, Pacific
                Northwest National Laboratory, Richland, WA 99352, USA}
\address[purdue1]{Department of Mathematics, Purdue University, West Lafayette, IN 47907, USA}
\address[purdue2]{School of Mechanical Engineering, Purdue University, West Lafayette, IN 47907, USA}
\cortext[gl]{Corresponding author.}

\ead{guanglin@purdue.edu}

\begin{abstract}
Compressive sensing has become a powerful addition to uncertainty quantification
in recent years. This paper identifies new bases for random variables through
linear mappings such that the representation of the quantity of interest is more
sparse with new basis functions associated with the new random variables. This
sparsity increases both the efficiency and accuracy of the compressive
sensing-based uncertainty quantification method. Specifically, we consider
rotation-based linear mappings which are determined iteratively for Hermite
polynomial expansions. We demonstrate the effectiveness of the new method with
applications in solving stochastic partial differential equations and
high-dimensional ($\mathcal{O}(100)$) problems.
\end{abstract}


\end{frontmatter}

$\dagger$ The first two authors made equal contributions to this manuscript.

\input{intro}

\input{review}

\input{method}

\input{numeric}

\input{conclusion}

\input{append}

\bibliographystyle{elsarticle-num}
\bibliography{uq}

\end{document}

%% file: mathhead.tex
\usepackage{amsmath}
\usepackage{amsfonts}
\usepackage{amssymb}
\usepackage{amsthm}
\usepackage{bm}
\usepackage{indentfirst}
\usepackage{titlesec}
\usepackage{graphicx}
\DeclareGraphicsExtensions{.eps}
\usepackage{subfigure}
\usepackage{array}
\usepackage{fullpage}
\usepackage{color}
\usepackage{algorithm}
\usepackage{algorithmic}
\usepackage{url}

\DeclareMathAlphabet{\mathsfsl}{OT1}{cmss}{m}{sl}
\newcommand{\tensor}[1]{\mathsfsl{#1}}

\newcommand{\dif}{\mathrm{d}}

\newcommand{\ba}{\bm\alpha}
\newcommand{\bx}{\bm\xi}
\newcommand{\ve}{\varepsilon}

\newcommand{\mexp}[1]{\mathbb{E}\left\{{#1}\right\}}
\newcommand{\mr}{\mathbb{R}}
\newcommand{\mn}{\mathbb{N}}

\newcommand{\PreserveBackslash}[1]{\let\temp=\\#1\let\\=\temp}
\newcolumntype{C}[1]{>{\PreserveBackslash\centering}p{#1}}
\newcolumntype{R}[1]{>{\PreserveBackslash\raggedleft}p{#1}}
\newcolumntype{L}[1]{>{\PreserveBackslash\raggedright}p{#1}}

\renewcommand{\theequation}{\thesection.\arabic{equation}}
\numberwithin{equation}{section}

\theoremstyle{definition}

\newcommand\Def{\stackrel{\textrm{def}}{=}}
\DeclareMathOperator{\sech}{sech}

%% file: intro.tex
\section{Introduction}
\label{sec:intro}
Uncertainty quantification (UQ) plays an important role in constructing
computational models as it helps to understand the influence of uncertainties on
the quantity of interest. In this paper, we study parametric uncertainty,
which treats some of the parameters as random variables.
Let $(\Omega,\mathcal{F},P)$ be a complete probability space, where $\Omega$ is
the event space and $P$ is a probability measure on the $\sigma$-field
$\mathcal{F}$. We consider a system depending on a $d$-dimensional random vector
$\bx(\omega)=(\xi_1(\omega),\xi_2(\omega),\cdots,\xi_d(\omega))^T$, where
$\omega$ is an event in $\Omega$. For simplicity, we denote $\xi_i(\omega)$ as
$\xi_i$. We aim to approximate the quantity of interest $u(\bx)$ with a 
generalized polynomial chaos (gPC) expansion \cite{GhanemS91, XiuK02}:
\begin{equation}\label{eq:gpc}
u(\bx) = \sum_{n=1}^Nc_n\psi_n(\bx) + \ve(\bx),
\end{equation}
where $\ve$ is the truncation error, $N$ is a positive integer, $c_n$ are
coefficients, $\psi_n$ are multivariate polynomials which are orthonormal with 
respect to the distribution of $\bx$:
\begin{equation}
\int_{\mr^d} \psi_i(\bx)\psi_j(\bx)\rho(\bx)\dif\bx = \delta_{ij},
\end{equation}
where $\rho(\bx)$ is the probability distribution function (PDF) of $\bx$ and 
$\delta_{ij}$ is the Kronecker delta. The approximation converges in the $L_2$ 
sense as $N$ increases if $u$ is in the Hilbert space associated with the 
measure of $\bx$ (i.e., the weight of the inner product is the PDF of $\bx$) 
\cite{XiuK02, CameronM1947, Ogura1972}. Stochastic Galerkin and probabilistic 
collocation are two popular methods 
\cite{GhanemS91,XiuK02,TatangPPM97, XiuH05,FooWK08,BabuskaNT10} used to 
approximate the gPC coefficients $\bm c=(c_1,c_2,\cdots,c_N)^T$. Stochastic
collocation  starts by generating samples of input $\bx^q, q=1,2,\cdots,M$
based on $\rho(\bx)$. Next, the computational model is calculated for each 
$\bx^q$ to obtain corresponding samples of the output $u^q=u(\bx^q)$. Finally,
$\bm c$ are approximated based on $u^q$ and $\bx^q$. Note that in many practical
problems, it is very costly to obtain $u^q$ and, due to the limited 
computational sources, we will often have $M<N$ or even $M\ll N$. The smaller
number of samples than basis functions implies that the following linear system
is under-determined:
\begin{equation}\label{eq:cs_eq}
\tensor\Psi \bm c = \bm u + \bm\ve,
\end{equation}
where $\bm u=(u^1,u^2,\cdots,u^M)^T$ is the vector of output samples, 
$\tensor\Psi$ is an $M\times N$ matrix with $\Psi_{ij}=\psi_j(\bx^i)$ and 
$\bm\ve=(\ve^1,\ve^2,\cdots,\ve^M)^T$ is a vector of error samples with 
$\ve^q=\ve(\bx^q)$. The compressive sensing method is effective at solving this
type of under-determined problem when $\bm c$ is sparse 
\cite{CandesRT06,DonohoET06,Candes08,BrucksteinDE09} and recent studies have 
applied this approach to uncertainty quantification (UQ) problems 
\cite{DoostanO11,YanGX12,YangK13,KaragiannisBL15,LeiYZLB14,PengHD14,XuZ14,
HamptonD15,TangI14,JakemanES14,SargsyanSNDRT14,PengHD15}. 

Several useful approaches have been developed to enhance the efficiency of 
solving Eq.~\eqref{eq:cs_eq} in UQ applications. First, re-weighted $\ell_1$ 
minimization assigns a weight to each $c_n$ and solves a weighted $\ell_1$ 
minimization problem to enhance the sparsity \cite{CandesWB08}. The weights can
be estimated in \textit{a priori} {\cite{PengHD14,RauhutW15}} or, for more general cases, 
can be obtained iteratively \cite{YangK13,LeiYZLB14}. Second, better sampling
strategies can be used, such as minimizing the mutual coherence 
{\cite{RauhutW12,HamptonD15}}. Third, Bayesian compressive sensing method provides the 
posterior distribution of the coefficients
\cite{SargsyanSNDRT14, KaragiannisBL15}. Finally, adaptive basis selection 
selects basis functions to enhance the efficiency instead of fixing the basis
functions at the beginning \cite{JakemanES14}. { Recently, we propose
an approach \cite{LeiYZLB14} to enhance the sparsity of $\bm c$ through the rotation of the random
vector $\bx$ to a new random vector $\bm\eta$, where the rotation operator is determined
by the sorted variability directions of the quantity of interest $u$ based on
the active subspace method \cite{ConstantineDW14}.}

In this work, we {aim to extend our previous work \cite{LeiYZLB14} and} consider the 
specific case where the system depends on 
i.i.d.~Gaussian random variables; i.e., $\bx\sim\mathcal{N}(\bm 0, \tensor I)$
where $\bm 0$ is a $d$-dimensional zero vector and $\tensor I$ is a $d\times d$
identity matrix. This assumption appears in a wide range of physics and 
engineering problems. We aim to find a mapping $g:\mr^d\mapsto\mr^d$ which maps
$\bx$ to a new set of i.i.d.~Gaussian random variables 
$\bm\eta=(\eta_1, \eta_2, \cdots, \eta_d)^T$ such that the gPC expansion of $u$
with respect to $\bm\eta$ is sparser. In other words,
\begin{equation}\label{eq:gpc2}
u(\bx) \approx \sum_{n=1}^N c_n\psi_n(\bx)=\sum_{n=1}^N\tilde
c_n\tilde\psi_n(\bm\eta(\bx))\approx u(\bm\eta(\bx)),
\end{equation}
where $\tilde\psi_n$ are orthonormal polynomials associated with the new random
vector $\bm\eta$ and $\tilde c_n$ are the corresponding coefficients. Note that
$\psi_n=\tilde\psi_n$ since $\bm\eta\sim\mathcal{N}(\bm 0, \tensor I)$. We 
intend to find the set 
$\tilde{\bm c}=(\tilde c_1, \tilde c_2, \cdots, \tilde c_N)^T$ which is 
sparser than $\bm c$ while preserving the properties of matrix 
$\tilde{\tensor\Psi}$ (with $\tilde{\Psi}_{ij}=\tilde\psi_j(\bm\eta^i)$)
close to those of $\tensor\Psi$ to improve the efficiency of the compressive 
sensing method. To accomplish this, we will use a linear mapping, based on the
idea of active subspaces \cite{ConstantineDW14}, to obtain $\bm\eta$ as first 
proposed in \cite{LeiYZLB14}. Unlike { our previous work}, we build this mapping 
iteratively in order to obtain a sparser $\tilde{\bm c}$ and improve the 
efficiency of the gPC approximation by compressive sensing. We also provide the
analytical form of the ``gradient matrix" (see Eq.\eqref{eq:grad_mat}) to avoid
estimating it with Monte Carlo methods. {Our method} is applicable for both
$\ell_0$ and $\ell_1$ minimization problems. Especially, for the latter, we can
also integrate the {present} method with re-weighted $\ell_1$ minimization method to 
further reduce the error. We demonstrate that, compared with the standard
compressive sensing methods, our approach reduces the relative $L_2$ error 
of the gPC approximation.

%% file: review.tex
\section{Brief review of the compressive sensing-based gPC method}
\label{sec:review}

\subsection{Hermite polynomial chaos expansions}
In this paper we study systems relying on $d$-dimensional Gaussian random
vector $\bx\sim\mathcal{N}(\bm 0, \tensor I)$. Therefore, the gPC basis
functions are constructed by tensor products of univariate orthonormal Hermite
polynomials. For a multi-index 
$\ba=(\alpha_1,\alpha_2,\cdots,\alpha_d), \alpha_i\in\mn\cup\{0\}$, we set
\begin{equation}\label{eq:tensor}
\psi_{\ba}(\bx) =
\psi_{\alpha_1}(\xi_1)\psi_{\alpha_2}(\xi_2)\cdots\psi_{\alpha_d}(\xi_d).
\end{equation}
For two different multi-indices 
$\ba_i=((\alpha_i)_{_1}, (\alpha_i)_{_2}, \cdots, (\alpha_i)_{_d})$ and 
$\ba_j=((\alpha_j)_{_1}, (\alpha_j)_{_2}, \cdots, (\alpha_j)_{_d})$, we have the
property
\begin{equation}
\int_{\mr^d} \psi_{\ba_i}(\bx)\psi_{\ba_j}(\bx) \rho(\bx) \dif\xi =
\delta_{\ba_i\ba_j} = \delta_{(\alpha_i)_{_1}(\alpha_j)_{_1}}
\delta_{(\alpha_i)_{_2}(\alpha_j)_{_2}}\cdots
\delta_{(\alpha_i)_{_d}(\alpha_j)_{_d}},
\end{equation}
where 
\begin{equation}
\rho(\bx) = \left(\dfrac{1}{\sqrt{2\pi}}\right)^{d}
\exp\left(-\dfrac{\xi_1^2+\xi_2^2+\cdots+\xi_d^2}{2}\right).
\end{equation}
For simplicity, we denote $\psi_{\ba_i}(\bx)$ as $\psi_i(\bx)$.


\subsection{Compressive sensing}

The vector $\bm c$ in Eq.~\eqref{eq:cs_eq} can be approximated by solving the 
following optimization problem:
\begin{equation}\label{eq:lh}
(P_{h,\epsilon}):~\arg \min_{\hat{\bm c}}\Vert\hat{\bm c}\Vert_h, 
\text{~~subject to~~} \Vert\tensor\Psi\hat{\bm c}-\bm u\Vert_2\leq\epsilon,
\end{equation}
where $\epsilon=\Vert\bm\ve\Vert_2$ and $h$ is typically set as $0$ or $1$. For
$h=0$ ($\ell_0$ minimization problem), the greedy  Orthogonal Matching Pursuit
(OMP) algorithm \cite{ChenDS98, BrucksteinDE09} can be applied; for $h=1$ 
($\ell_1$ minimization problem), convex optimization methods are directly
applicable \cite{BoydV04}. As pointed out in \cite{BrucksteinDE09}, OMP is very
efficient -- when it works -- but convergence to a sparse solution is not always
guaranteed. There are specific cases where a sparse solution is possible while 
OMP yields a dense one. Since both the OMP and $\ell_1$ minimization approaches
are widely used, we will demonstrate the effectiveness of our new method for 
both methods.

Next, we introduce the concept of \emph{sparsity} as it is critical in the error
estimates for solving the under-determined system Eq.~\eqref{eq:cs_eq} with the
compressive sensing method. The $\ell_0$ ``norm" of vector
$\bm x=(x_1,x_2,\cdots,x_N)$ is defined as the number of its non-zeros entries 
\cite{Donoho06,CandesRT06,BrucksteinDE09}
\begin{equation}
\Vert\bm x\Vert_0\Def \#\{i:x_i\neq 0\}
\end{equation}
and $\ell_1$ norm is defined as the sum of the absolute value of its entries:
\begin{equation}
\Vert\bm x\Vert_1\Def \sum_{n=1}^N |x_n|.
\end{equation}
$\bm x$ is called \emph{$s$-sparse} if $\Vert \bm x\Vert_0\leq s$, and $\bm x$ 
is considered a sparse vector if $s\ll N$. Few practical systems have a truly 
sparse gPC coefficients $\bm c$. However, in many cases, the $\bm c$ are 
compressible, i.e., only a few entries make significant contribution to its
$\ell_1$ norm. In subsequent discussion, we relax the definition of ``sparse":
$\bm x$ is considered sparse if $\Vert \bm x - \bm x_s\Vert_1$ is small for 
$s\ll N$. Here $\bm x_s$ is defined as the 
\textbf{best $s$-sparse approximation} one could obtain if one knew exactly the
locations and amplitudes of the $s$-largest entries of $\bm x$, i.e., $\bm x_s$
is the vector $\bm x$ with all but the $s$-largest entries set to 
zero \cite{Candes08}.

The error bound for solving Eq.~\eqref{eq:cs_eq} with $\ell_1$ minimization 
requires definition of the \textit{restricted isometry property} (RIP) constant
\cite{CandesT05}. For each integer $s=1,2,\cdots$, the isometry constant 
$\delta_s$ of a matrix $\tensor\Phi$ is defined as the smallest number such that
\begin{equation}
(1-\delta_s)\Vert\bm x\Vert_2^2\leq \Vert\tensor\Phi\bm x\Vert_2^2\leq
(1+\delta_s)\Vert\bm x\Vert_2^2
\end{equation}
holds for all $s$-sparse vectors $\bm x$.
\noindent With some restrictions, Candes et al.\ showed $\bm x$ can be stably
reconstructed \cite{Candes08}.
Assume that the matrix $\tensor\Psi$ satisfies $\delta_{2s}<\sqrt{2}-1$, and
$\Vert\bm\ve\Vert_2\leq\epsilon$, then solution 
$\hat{\bm c}$ to $(P_{1,\epsilon})$ obeys
\begin{equation}\label{eq:l1_thm}
\Vert \bm c - \hat{\bm c}\Vert_2 \leq C_1\epsilon + C_2 
\dfrac{\Vert\bm c-\bm c_s\Vert_1}{\sqrt{s}},
\end{equation}
where $C_1$ and $C_2$ are constants, $\bm c$ is the exact vector we aim to 
approximate and $\hat{\bm c}$ is the solution of $(P_{1,\epsilon})$. 
This result implies that the upper bound of the error is related to the 
truncation error and the sparsity of $\bm c$, which is indicated in the first
and second terms on the right hand side of Eq. \eqref{eq:l1_thm}, respectively.

The \emph{re-weighted} $\ell_1$ minimization approach is an improvement of the
$\ell_1$ minimization method, which enhances the accuracy of estimating $\bm c$
\cite{CandesWB08}. The re-weighted $\ell_1$ approach solves the following 
optimization problem:
\begin{equation}\label{eq:wl1}
(P_{1,\epsilon}^W):~\arg\min_{\hat{\bm c}}\Vert\tensor W\hat{\bm c}\Vert_1, 
~\text{subject to}~\Vert\tensor\Psi\hat{\bm c}-\bm u\Vert_2\leq\epsilon,
\end{equation}
where $\tensor W$ is a diagonal matrix: 
$\tensor W=\text{diag}(w_1,w_2,\cdots,w_N)$. Clearly, $(P_{1,\epsilon})$ can be
considered as a special case of $(P_{1,\epsilon}^W)$ by setting 
$\tensor W=\tensor I$. The elements $w_i$ of the diagonal matrix can be 
estimated based on analysis of $u$ as in Peng et al.~\cite{PengHD14}, or be
estimated iteratively \cite{CandesWB08, YangK13}. More precisely, for each 
iteration $l$, $(P_{1,\epsilon}^W)$ is solved to obtain $\hat{\bm c}^{(l)}$ and
then $w_i^{(l+1)}=1/(|\hat c_i^{(l)}|+\delta)$ for the next step. The parameter 
$\delta>0$ is introduced to provide stability and to ensure that a zero-valued
component in $\hat{\bm c}^{(l)}$ does not prohibit a nonzero estimate at the
next step. In Candes et al.~\cite{CandesWB08}, the authors suggest two to three
iterations of this procedure. Subsequent analytical work \cite{Needell09} 
provides an error bound for each iteration as well as the limit of computing
$\hat{\bm c}$ with re-weighted $\ell_1$ minimization. The form is similar to
Eq.~\eqref{eq:l1_thm} with different constants.

In practice, the error term $\epsilon$ is not known \emph{a priori}, hence
cross-validation is needed to estimate it. One such algorithm is 
\cite{DoostanO11} summarized in Algorithm \ref{algo:cross} :
\begin{algorithm}
\caption{Cross-validation to estimate the error $\epsilon$}
\label{algo:cross}
\begin{algorithmic}[1]
\STATE Divide the $M$ output samples to $M_r$ reconstruction ($\bm u_r$) and
$M_v$ validation ($\bm u_v$) samples and divide the measurement matrix 
$\tensor\Psi$ correspondingly into $\tensor\Psi_r$ and $\tensor\Psi_v$.
\STATE Choose multiple values for $\epsilon_r$ such that the exact error
$\Vert\tensor\Psi_r\bm c-\bm u_r\Vert_2$ of the reconstruction samples is
within the range of $\epsilon_r$ values.
\STATE For each $\epsilon_r$, solve $(P_{h,\epsilon})$ with $\bm u_r$ and
$\tensor\Psi_r$ to obtain $\hat{\bm c}$, then compute 
$\epsilon_v=\Vert\tensor\Psi_v\hat{\bm c}-\bm u_v\Vert_2$.
\STATE Find the minimum value of $\epsilon_v$ and its corresponding 
$\epsilon_r$. Set $\epsilon=\sqrt{M/M_r}\epsilon_r$.
\end{algorithmic}
\end{algorithm}

We omit the review of the theoretical results for the OMP as well as its 
variants, and refer interested readers to the literature 
\cite{BrucksteinDE09,Tropp04,Tropp07}. Similar to the $\ell_1$ approach, the 
error estimate for OMP includes a term which depends on the sparsity of $\bm c$.
This is a critical point that motivates us to propose the new method described
in the next section.


\subsection{Compressive sensing-based gPC methods}
Given $M$ samples of $\bx$, the quantity of interest $u$ is approximated by a 
gPC expansion as in Eq.~\eqref{eq:gpc}:
\begin{equation}
u(\bx^q) = \sum_{n=1}^N c_n\psi(\bx^q) + \ve(\bx^q), \quad q=1,2,\cdots,M,
\end{equation}
which can be rewritten as Eq.~\eqref{eq:cs_eq}. A typical approach to compressive
sensing based-gPC is summarized in Algorithm \ref{algo:cs1}.
\begin{algorithm}
\caption{Compressive sensing-based gPC}
\label{algo:cs1}
\begin{algorithmic}[1]
\STATE Generate input samples $\bx^q, q=1,2,\cdots, M$ based on the distribution
of $\bx$. 
\STATE Generate output samples $u^q=u(\bx^q)$ by solving the complete model; 
e.g., running simulations, solvers, etc.
\STATE Select gPC basis functions $\{\psi_n\}_{n=1}^N$ associated with $\bx$ and
then generate the measurement matrix $\tensor\Psi$ by setting
$\Psi_{ij}=\psi_j(\bx^i)$.
\STATE Solve the optimization problem $(P_{h,\epsilon})$:
\[\arg \min_{\hat{\bm c}}\Vert\hat{\bm c}\Vert_h, ~
\text{subject to} \Vert\tensor\Psi\hat{\bm c}-\bm u\Vert_2\leq\epsilon,\]
where $h=0$ or $1$, $\bm u=(u^1,u^2,\cdots,u^M)^T$, and $\epsilon$ is obtained by
cross-validation. If the re-weighted $\ell_1$ method is employed, solve
$(P_{1,\epsilon}^W)$ instead.
\STATE Set $\bm c=\hat{\bm c}$ and construct gPC expansion as 
$u(\bx)\approx \sum_{n=1}^N c_n\psi_n(\bx)$.
\end{algorithmic}
\end{algorithm}

Note that the RIP condition in Theorem 2.2 is sufficient but not necessary; 
furthermore, it is difficult to obtain the exact RIP constant in practical 
problems. A more tractable property of the measurement matrix for calculation is
the \emph{mutual coherence} \cite{BrucksteinDE09}:
\begin{equation}
\mu(\tensor\Psi) = \max_{1\leq j,k\leq N,j\neq k}
\dfrac{|\bm\Psi_j^T\bm\Psi_k|}{\Vert\bm\Psi_j\Vert_2\cdot\Vert\bm\Psi_k\Vert_2},
\end{equation}
where $\bm\Psi_j$ and $\bm\Psi_k$ are columns of $\tensor\Psi$. In general, a
measurement matrix with smaller mutual coherence is better able to recover a 
sparse solution with the compressive sensing method. Note that
$\mexp{\psi_i(\xi)\psi_j(\xi)}=\delta_{ij}$ since $\{\psi_i\}_{i=1}^N$ are 
orthonormal polynomials. Therefore, asymptotically, $\mu(\tensor\Psi)$ converges
to zeros according to the strong law of large numbers. 

In the next section, we will demonstrate that our new method increases the 
sparsity of $\bm c$ without changing $\mu$ significantly, { and} hence, our method
is able to improve the accuracy of the compressive sensing-based gPC method.

%% file: method.tex
\section{Iterative rotations for increasing sparsity}
\label{sec:method}

In this section, we provide a heuristic method to identify the rotation matrix
by computing the eigenvalue decomposition of a gradient matrix $\tensor G$. 
The rotation increases the sparsity of the gPC expansions of quantity of 
interest $u$ with respect to a new set of random variables. The rotation 
procedure is applied iteratively to achieve a target sparsity level. The
enhancement of the sparsity decreases the second term (sparsity-induced error)
on the right hand side of Eq.~\eqref{eq:l1_thm}. For the cases where this 
sparsity-induced error dominates the total error of the compressive sensing 
method, our new approach improves the overall accuracy. 

From Eq.~\eqref{eq:l1_thm}, we notice that if $\bm c$ is exactly sparse (i.e.,
$\bm c=\bm c_{s^*}$ for some $s^*\ll N$) and if the RIP condition is satisfied 
for $s\geq s^*$, then $\Vert \bm c-\bm c_s\Vert=0$. Therefore, the upper bound
of the error only depends on $\epsilon$. In practical problems, $\bm c$ is 
usually not exactly sparse. But if the truncation error {$\epsilon$} is 
sufficiently small, then the second term on the right hand side of 
Eq.~\eqref{eq:l1_thm} dominates the upper bound of the error. Hence, in order 
to improve the accuracy of the gPC expansion, we need to decrease 
$\Vert\bm c-\bm c_s\Vert_1/\sqrt{s}$. However, once the gPC basis functions are
selected, $\bm c$, and therefore $\Vert\bm c-\bm c_s\Vert_1/\sqrt{s}$, are 
fixed. A natural way to enhance the sparsity of $\bm c$ is to find another set 
of random variables $\bm\eta=(\eta_1,\eta_2,\cdots,\eta_{\tilde d})^T$, which 
depend on $\bx$ such that the vector $\tilde{\bm c}$, which are the gPC 
coefficients of $u$ with respect to $\bm\eta$, is sparser. In order words, our 
goal is to seek $\bm\eta(\bx)$ with
\begin{equation*}
u(\bx) \approx \sum_{n=1}^Nc_n\psi_n(\bx)=\sum_{n=1}^{\tilde N}\tilde
c_{n}\tilde\psi_{n}(\bm\eta(\bx))\approx u(\bm\eta(\bx)),
\end{equation*}
such that $\Vert\tilde{\bm c}-\tilde{\bm c}_s\Vert_1<\Vert\bm c-\bm c_s\Vert_1$. 
Note that $\tilde N$ does not necessarily equal $N$ and $\tilde d$ can be
different from $d$. We will denote the mapping from $\bx$ to
$\bm\eta$ as $g:\mr^{d}\mapsto\mr^{\tilde d}$.

There are several behaviors { that} our proposed approach must exhibit.
\begin{itemize}
\item \emph{The PDF of $\bm\eta$ must be computed efficiently.} The first step
of generating a new gPC expansion is to obtain the PDF of $\bm\eta$. Hence, if
$g$ is complicated, the PDF of $\bm\eta$ will be difficult to obtain. Even if
$g$ is a simple function, it can still be difficult to obtain an accurate PDF if
the dimension is large.
\item \emph{The new gPC basis functions associated with $\bm\eta$ must be 
computed efficiently.} If the $\eta_i$ are independent, then the new gPC basis
can be constructed as the tensor product of univariate basis functions in each
dimension. Although this is not necessary, it will make the construction of new
basis functions easier as it avoids the computation of high-dimensional 
integrals. 
\item \emph{The properties of the measurement matrix must be preserved.} 
Clearly, the measurement matrix changes as we introduce new random variables and
new basis functions. Even though we may construct a very sparse $\tilde{\bm c}$,
if the key properties of the measurement matrix are altered too much (e.g., the 
RIP constant or mutual coherence increases dramatically), we may be unable to 
obtain an accurate result with the compressive sensing method.
\item \emph{No additional output samples must be needed.} In particular, the 
existing output samples $u^q,q=1,2,\cdots, M$ should be sufficient. This is 
especially important for the cases when the model (simulation or deterministic 
solver) is very costly to compute.
\end{itemize}

In this work, we focus on the special case of normal distributions
$\bx\sim\mathcal{N}(\bm 0,\tensor I)$; hence, the $\psi_i$ are constructed as 
the tensor product of univariate orthonormal Hermite polynomials as shown in
Eq.~\eqref{eq:tensor}. We aim to find a linear mapping $g:\mr^d\mapsto\mr^d$ 
such that 
\begin{equation}\label{eq:linear_map}
\bm\eta=g(\bx)=\tensor A\tensor\bx,
\end{equation}
where $\tensor A$ is an orthonormal matrix.

If we find this { matrix $\tensor A$}
then all of the aforementioned behaviors can be 
obtained. We know that $\bm\eta\sim\mathcal{N}(\bm 0, \tensor I)$ since 
$\tensor A\tensor A^T=\tensor I$. Therefore, the samples of $\bm\eta$ can be
obtained as $\bm\eta^q=\tensor A\bx^q$, where $\bx^q$ are generated at the 
beginning (Step 1 in Algorithm \ref{algo:cs1}). Since 
$\bm\eta\sim\mathcal{N}(\bm 0,\tensor I)$ we can set $\tilde\psi_i=\psi_i$ and
no additional computation is needed. The difference between $\tensor\Psi$ and
$\tilde{\tensor\Psi}$ is that the latter is constructed by evaluating 
orthonormal Hermite polynomials at another set of samples of i.i.d.~Gaussian 
random variables; i.e., 
$\tilde\Psi_{ij}=\tilde\psi_j(\bm\eta^i)=\psi_j(\bm\eta^i)$. Therefore, the
mutual coherence of $\tilde{\tensor\Psi}$ converges to 0 as that of
$\tensor\Psi$, and the difference between $\mu(\tensor\Psi)$ and
$\mu(\tilde{\tensor\Psi})$ is $\mathcal{O}(M^{-1/2})$, the deviation of the 
Monte Carlo numerical integral from the exact value. No additional samples 
$u^q$ are required since the improvement of accuracy is achieved by enhancing
the sparsity of gPC coefficients.

Given the Hermite polynomials defined above, we have a new expansion for $u$:
\begin{equation}
u(\bx) \approx \sum_{n=1}^Nc_n\psi_n(\bx)=
\sum_{n=1}^N\tilde c_n\psi_n(\tensor A\bx))\approx u(\bm\eta)
\end{equation}
with $\tilde{\bm c}$ sparser than $\bm c$.
In order to obtain the $\tensor A$, we adopt the active subspace approach
\cite{ConstantineDW14}. We first define the ``gradient matrix":
\begin{equation}\label{eq:grad_mat}
\tensor G \Def
\mexp{\nabla u(\bx)\cdot\nabla u(\bx)^T} = \tensor U\tensor \Lambda\tensor U^T,
  \quad \tensor U\tensor U^T = \tensor I,
\end{equation}
where $\tensor G$ is symmetric, 
$\nabla u(\bx) = (\partial u/\partial \xi_1, \partial u/\partial \xi_2,
\cdots, \partial u/\partial \xi_d)^T$ is a column vector, 
$\tensor U=(\bm U_1,\bm U_2,\cdots, \bm U_d)$ is an orthonormal matrix
consisting of eigenvectors $\bm U_i$, and 
$\tensor\Lambda=\text{diag}(\lambda_1,\lambda_2,\cdots,\lambda_d)$ with 
$\lambda_1\geq\lambda_2\geq\cdots$ is a diagonal matrix with elements 
representing decreasing variation of the system along the respective 
eigenvectors. We choose $\tensor A=\tensor U^T$ which, as a unitary matrix,
defines a rotation in $\mr^d$ {and the linear mapping $g$ in
Eq.~\eqref{eq:linear_map} projects $\bx$ on the eigenvectors $\bm U_i$ }. 
Consequently, when the differences between { $|\lambda_i|$} are large, $g$ helps to
concentrate the dependence of $u$ primarily on the first few new random 
variables $\eta_i$ due to the larger variation of $u$ along the directions of
the corresponding eigenvectors. Therefore, we obtain a sparser $\tilde{\bm c}$
than $\bm c$. We note that this approach of constructing $\tensor G$ is similar
to the method of outer product gradients (OPGs) in statistics 
\cite{HardinH07, Xia07}. The information of the gradient of $u$ is also utilized
to improve the efficiency of compressive sensing in the
\emph{gradient-enhanced} method
\cite{RoderickAF10,LiARH11,LockwoodM13,JakemanES14,PengHD15}.

Since $u$ is not known \emph{a priori}, we replace it with its gPC expansion 
$u_{g}=\sum_{n=1}^Nc_n\psi_n(\bm\xi)$. In prior work by Constantine and others
\cite{ConstantineDW14, LeiYZLB14}, the expectation is obtained by taking the
average of the Monte Carlo results. In the current work, we compute $\tensor G$
differently: after obtaining $\bm c$ with compressive sensing method, we 
construct a gPC approximation $u_g$ to $u$ and approximate $\tensor G$ 
accordingly:
\begin{equation}
\tensor G\approx\mexp{\nabla\left(\sum_{n=1}^Nc_n\psi_n(\bm\xi)\right) 
    \cdot \nabla\left(\sum_{n'=1}^Nc_{n'}\psi_{n'}(\bm\xi)\right)^T}.
\end{equation}
The entries of $\tensor G$ can be approximated as:
\begin{equation}\label{eq:grad}
\begin{aligned}
G_{ij} & \approx\mexp{
\dfrac{\partial}{\partial\xi_i}\left(\sum_{n=1}^N c_n\psi_n(\bm\xi)\right)\cdot
\dfrac{\partial}{\partial\xi_j}\left(\sum_{n'=1}^Nc_{n'}\psi_{n'}(\bm\xi)\right)}
\\ & =
\mexp{\left(\sum_{n=1}^Nc_n\dfrac{\partial\psi_n(\bm\xi)}{\partial\xi_i}\right)
\cdot
\left(\sum_{n'=1}^Nc_{n'}\dfrac{\partial\psi_{n'}(\bm\xi)}{\partial\xi_j}\right)}
\\ & =
\sum_{n=1}^N\sum_{n'=1}^Nc_nc_{n'}
\mexp{\dfrac{\partial\psi_n(\bm\xi)}{\partial\xi_i}\cdot
\dfrac{\partial\psi_{n'}(\bm\xi)}{\partial\xi_j}} \\
& = \bm c^T \tensor K_{ij} \bm c,
\end{aligned}
\end{equation}
where $\tensor K_{ij}$ is a ``stiffness'' matrix with entries
\begin{equation}\label{eq:kernel}
(K_{ij})_{kl} = \mexp{\dfrac{\partial\psi_k(\bm\xi)}{\partial\xi_i}
                \cdot\dfrac{\partial\psi_{l}(\bm\xi)}{\partial\xi_j}}. 
\end{equation}                
Notice that $\tensor K_{ij}$ can be precomputed since $\{\psi_i\}$ are normalized
Hermite polynomials (see Appendix for details). $\tensor G$ is a $d\times d$
matrix, where $d$ is the number of random variables in the system. Since 
$\tensor G$ is a symmetric matrix, we only need to compute $d(d+1)/2$ of its 
entries. Furthermore, unlike the active subspace method, which focuses on the 
\emph{subspace} of $\mr^d$, we keep the dimension and set of basis functions
unchanged.

The entire iterative procedure is summarized in Algorithm \ref{algo:cs_rot}. 
\begin{algorithm}
\caption{Compressive sensing method with iterative rotations}
\label{algo:cs_rot}
\begin{algorithmic}[1]
\STATE For given random vector $\bx$ and quantity of interest $u$, run Algorithm
\ref{algo:cs1} to obtain approximated gPC coefficients $\hat{\bm c}$.
\STATE Set counter $l=0$, $\eta^{(0)}=\bx$, $\tilde{\bm c}^{(0)}=\hat{\bm c}$.
\STATE Construct $\tensor G^{l+1}$ with $\hat{\bm c}^{(l)}$ according to 
Eq.~\eqref{eq:grad}. Then decompose $\tensor G^{(l+1)}$ as
\[\tensor G^{(l+1)}=
\tensor U^{(l+1)}\tensor\Lambda^{(l+1)}(\tensor U^{(l+1)})^T,\quad
\tensor U^{(l+1)}(\tensor U^{(l+1)})^T=\tensor I.\]
\STATE Define $\bm\eta^{(l+1)}=(\tensor U^{(l+1)})^T\bm\eta^{(l)}$, and compute
samples 
$(\bm\eta^{(l+1)})^q=(\tensor U^{(l+1)})^T(\bm\eta^{(l)})^q, q=1,2,\cdots,M$.
Also, construct the new measurement matrix $\tensor\Psi^{(l+1)}$ with
$\Psi^{(l+1)}_{ij}=\psi_j((\bm\eta^{(l+1)})^i)$. 
\STATE Solve the optimization problem $(P_{h,\epsilon^{(l+1)}})$:
\[\arg \min_{\hat{\bm c}}\Vert\hat{\bm c}\Vert_h, \quad\text{subject to}
\Vert\tensor\Psi^{(l+1)}\hat{\bm c}-\bm u\Vert_2\leq\epsilon^{(l+1)},\]
and set $\tilde{\bm c}^{(l+1)}=\hat{\bm c}$. If reweight $\ell_1$ method is
employed, solve $(P_{1,\epsilon^{(l+1)}}^W)$ instead.
\STATE Set $l=l+1$. 
If $\left|\Vert\tensor U^{(l)}\Vert_1-d\right|<\theta$, where the threshold
$\theta$ is a positive real number, then stop. Otherwise, go to Step 3.
\end{algorithmic}
\end{algorithm}
This algorithm adds post-processing steps to Algorithm \ref{algo:cs1}, which
{ is} 
designed to increase the accuracy of compressive sensing based gPC method. In
Step 5, we use notation $\epsilon^{(l+1)}$ since the estimated error at 
iteration $l+1$ may be different from $\epsilon$. According to our numerical 
experiments (see Sec.~\ref{sec:numeric}), it is usually sufficient to test two
or three different values on $[\epsilon/5,\epsilon]$ in the cross-validation 
procedure (see Algorithm \ref{algo:cross}) to obtain $\epsilon^{(l+1)}$. 

In Algorithm \ref{algo:cs_rot}, we propose a terminating condition based on the
$\ell_1$ norm of the rotation matrix in each iteration:
$S(\tensor U^{(l)})\Def\sum_{i=1}^d\Vert\bm U^{(l)}_i\Vert_1$. If 
$\tensor U^{(l)}$ is the identity matrix or a permutation matrix, we need no
further iterations, and $S(\tensor U^{(l)})=d$. Otherwise, 
$S(\tensor U^{(l)})>d$ since $\Vert\bm U^{(l)}_i\Vert_2=1$ and
\begin{equation}
\Big\Vert\bm U^{(l)}_i\Big\Vert_1^2
=\left(\sum_{j=1}^d \left|(U_i^{(l)})_j\right|\right)^2 
=\Big\Vert\bm U_i^{(l)}\Big\Vert_2^2 + 2\sum_{1\leq j,k\leq d, j\neq k}
\left|(U_i^{(l)})_j\right|\left|(U_i^{(l)})_k\right|>1.
\end{equation}
Hence, one may set a threshold $\theta$ and the iteration stops when
$|S(\tensor U^{(l)})-d|<\theta$. Empirically, $\theta$ can be set around
$0.1d\sim 0.2d$. More sophisticated terminating conditions (e.g., sparsity or 
$\epsilon$ estimates) are also possible. 
{ A rigorous theoretical analysis on the convergence 
behavior is not available at this time. The criterion presented here provides an
approach to estimate, to some extend, whether our method converges. Empirically,
when this stopping criterion is satisfied, additional iterations will not
improve the accuracy significantly.} 
We {
also} note that the simplest terminating condition in Step 6 is to set a maximum 
iteration steps $L$. Based on our numerical examples in Sec.~\ref{sec:numeric}, 
this simple condition can also be useful. { In general, the 
efficiency of our method depends on the intrinsic sparsity of the system, i.e., 
whether the system mainly relies on a small amount of subspaces. The fewer
subspaces the system depends on, the better performance our method exhibits.
Otherwise, this method is less effective, e.g., an extreme case is
$u(\bm\xi)=\sum_{i=1}^d\xi_i^2$, for which the iterative rotations based on
current framework does not help.
}

%% file: numeric.tex
\section{Numerical results}
\label{sec:numeric}

In this section, we present five numerical examples to demonstrate the
effectiveness of our new method. The accuracies of different methods are
measured by the relative $L_2$ error:
$(\Vert u - u_g\Vert_2)/\Vert u\Vert_2$, where $u_g$ is the Hermite
polynomial expansion of $u$. The integral
\begin{equation}
\Vert u(\bx)\Vert_2 = \left(\int_{\mr^d} u(\bx)^2\rho(\bx)\dif\bx \right)^{1/2}
\end{equation}
(and $\Vert u-u_g\Vert_2$) is approximated with a high-level sparse grid method 
which is based on one-dimensional Gauss-Hermite quadrature and the Smolyak 
structure \cite{Smolyak63}. The term ``level'' $p$ means that the algebraic 
accuracy of the sparse grid method is $2p-1$. We use $P$ to denote the 
truncation order, which implies that Hermite polynomials up to order $P$
{ are}
included in expansion $u_g$. Hence, the number of unknowns can be computed as
$N=\bigl(\begin{smallmatrix}P+d\\d\end{smallmatrix}\bigr)$.

The relative errors we present in this section are obtained from $100$ 
independent replicates for each sample size $M$. For example, we generate $100$
independent sets of input samples $\bx^q,q=1,2,\cdots,M$, compute $100$ 
different relative errors, and then report the average of these error samples.
To investigate the effectiveness of the increasing of output samples, we set 
the $x$-axis in our figures as the ratio $M/N$ which is the fraction of 
available data with respect to number of unknowns. We use MATLAB package 
\texttt{SPGL1} \cite{BergF08,spgl1} to solve $(P_{1,\epsilon})$ as well as
$(P_{1,\epsilon}^W)$ and use \texttt{SparseLab} \cite{DonohoDST} for the OMP
method. If not otherwise indicated, results are obtained with $L=3$ iterations
in Step 6 of Algorithm \ref{algo:cs_rot}.

\subsection{Example function with equally important random variables} 
\label{subsec:ex1}
Consider the following function
\begin{equation}\label{eq:ex1}
u(\bx) = \sum_{i=1}^d \xi_i + 0.25\left(\sum_{i=1}^d \xi_i\right)^2
       + 0.025\left(\sum_{i=1}^d \xi_i\right)^3,
\end{equation}
where all $\xi_i$ are equally important. In this case, adaptive methods that 
build the surrogate model hierarchically based on the importance of $\xi_i$ 
(e.g., \cite{MaZ09,MaZ10,YangCLK12, ZhangYMMEKD14}) may not be efficient. 
A simple rotation matrix for this example has the form
\[\tensor A = 
\begin{pmatrix}
d^{-1/2} & d^{-1/2} & \cdots & d^{-1/2} \\
 &   &  & \\
 & \multicolumn{2}{c}{\tilde{\tensor{A}}}  & \\
 &   &  & 
\end{pmatrix},\]
where $\tilde{\tensor A}$ is a $(d-1)\times d$ matrix chosen to ensure that
$\tensor A$ is orthonormal. Given this choice for $\tensor A$, then
$\eta_1=(\sum_{i=1}^d\xi_i)/d^{1/2}$ and $u$ has a very simple representation:
\[u(\bx)=u(\bm\eta)=d^{1/2}\eta_1+0.25d\eta_1^2+0.025d^{3/2}\eta_1^3.\]
Therefore, as we keep the set of the basis functions unchanged, all the Hermite
polynomials not related to $\eta_1$ make no contribution to the expansion, which
implies that we obtain a very sparse representation of $u$. Unfortunately, the 
optimal structure is not known \textit{a priori}, hence, the standard 
compressive sensing cannot take advantage of it.

\begin{figure}[t]
\centering
\includegraphics[width=3in]{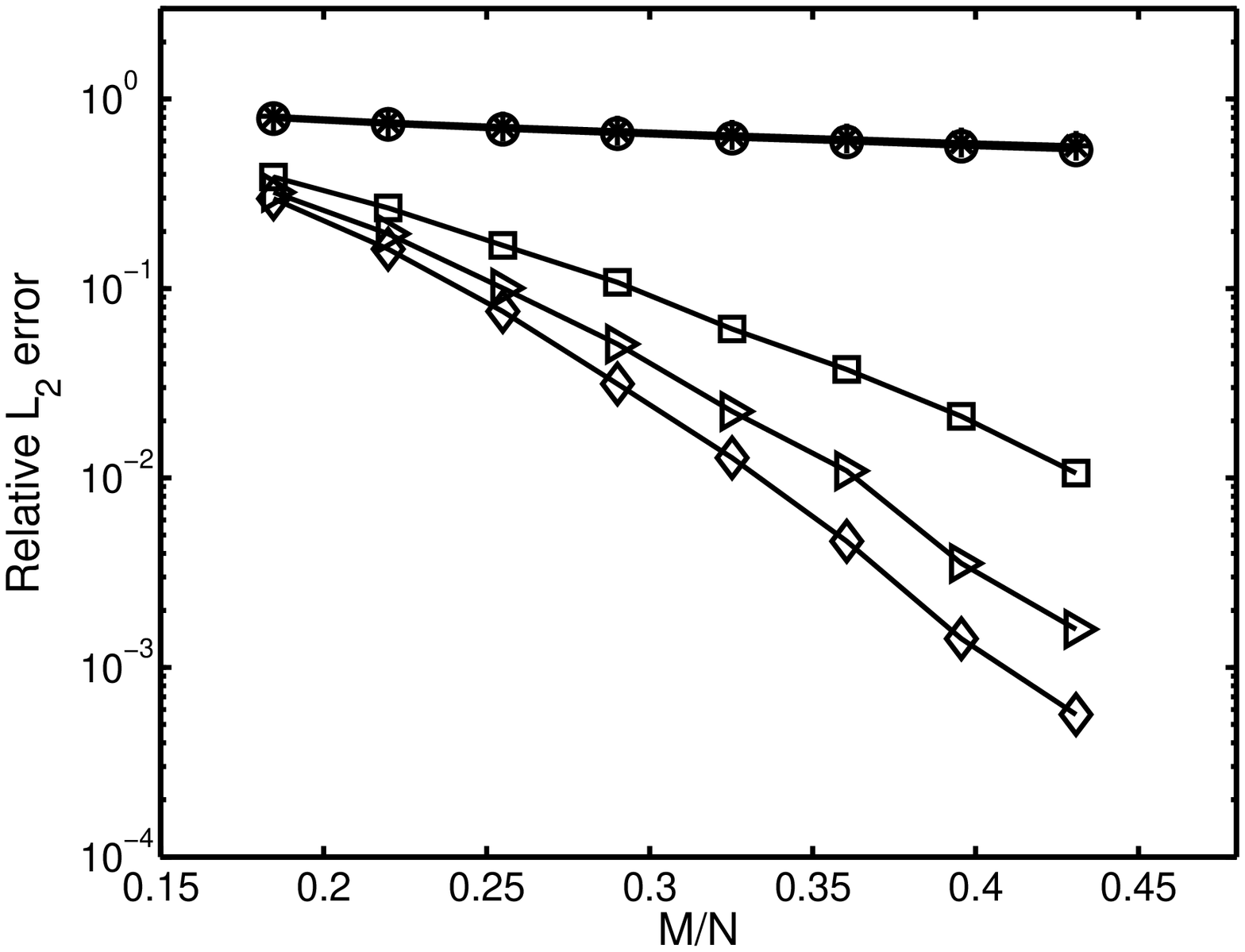}
\includegraphics[width=3in]{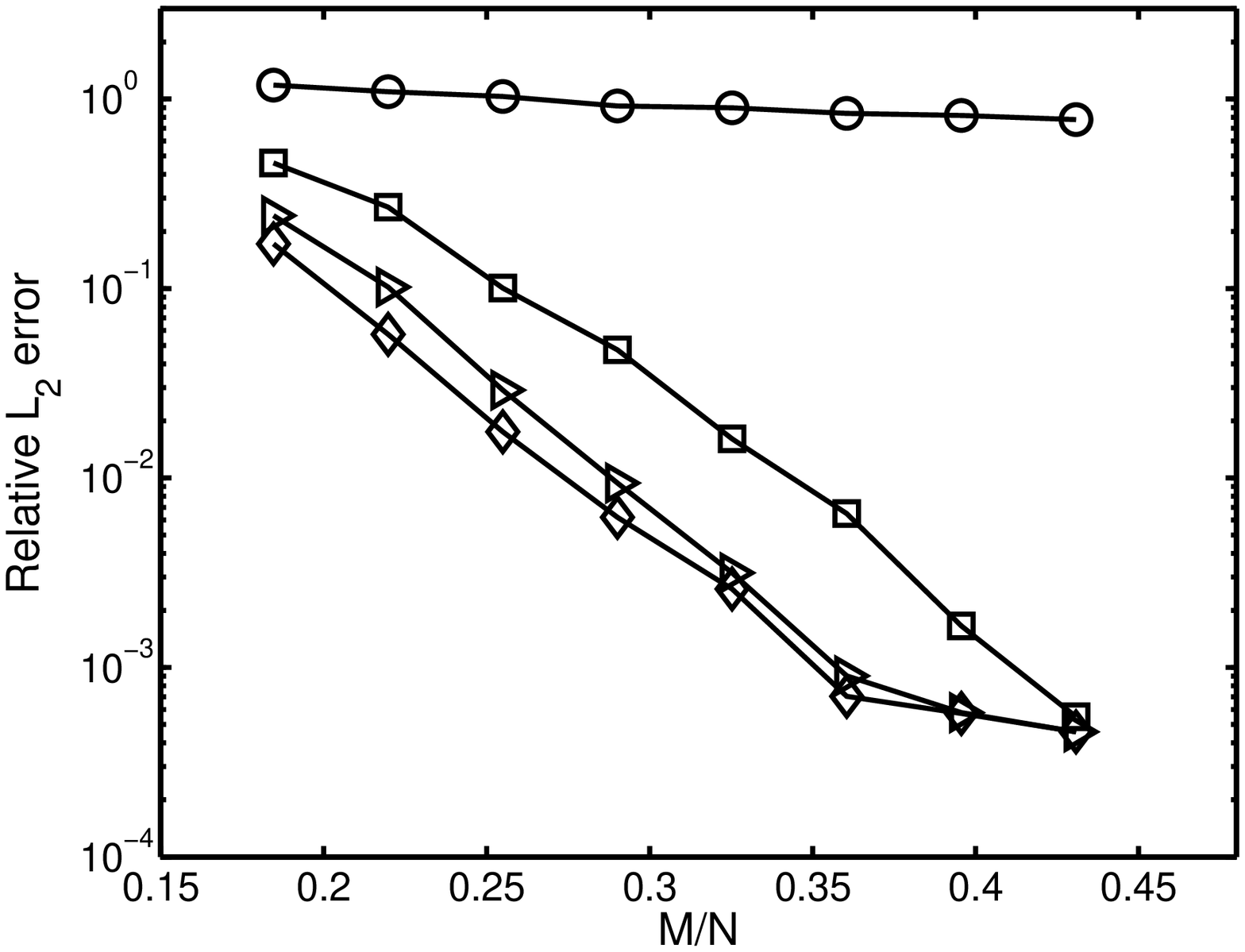}
\caption{Results for the example function with equally important random 
variables. (Left) Comparison with the $\ell_1$ method.  ``$\circ$": $\ell_1$,
``$\Box$": rotated $\ell_1$ with $3$ iterations, ``$\triangleright$": rotated
$\ell_1$ with $6$ iterations, ``$\diamond$": rotated $\ell_1$ with $9$ 
iterations, ``$*$": re-weighted $\ell_1$. (Right) Comparison with the OMP 
method. ``$\circ$": OMP, ``$\Box$": rotated OMP with $3$ iterations, 
``$\triangleright$": rotated OMP with $6$ iterations, ``$\diamond$": rotated 
OMP with $9$ iterations. These calculations were performed with dimension 
$d=12$ and the number of unknowns $N=455$.}
\label{fig:ex1_l1} \label{fig:ex1_omp}
\end{figure}

\begin{figure}[t]
\centering
\includegraphics[width=3in]{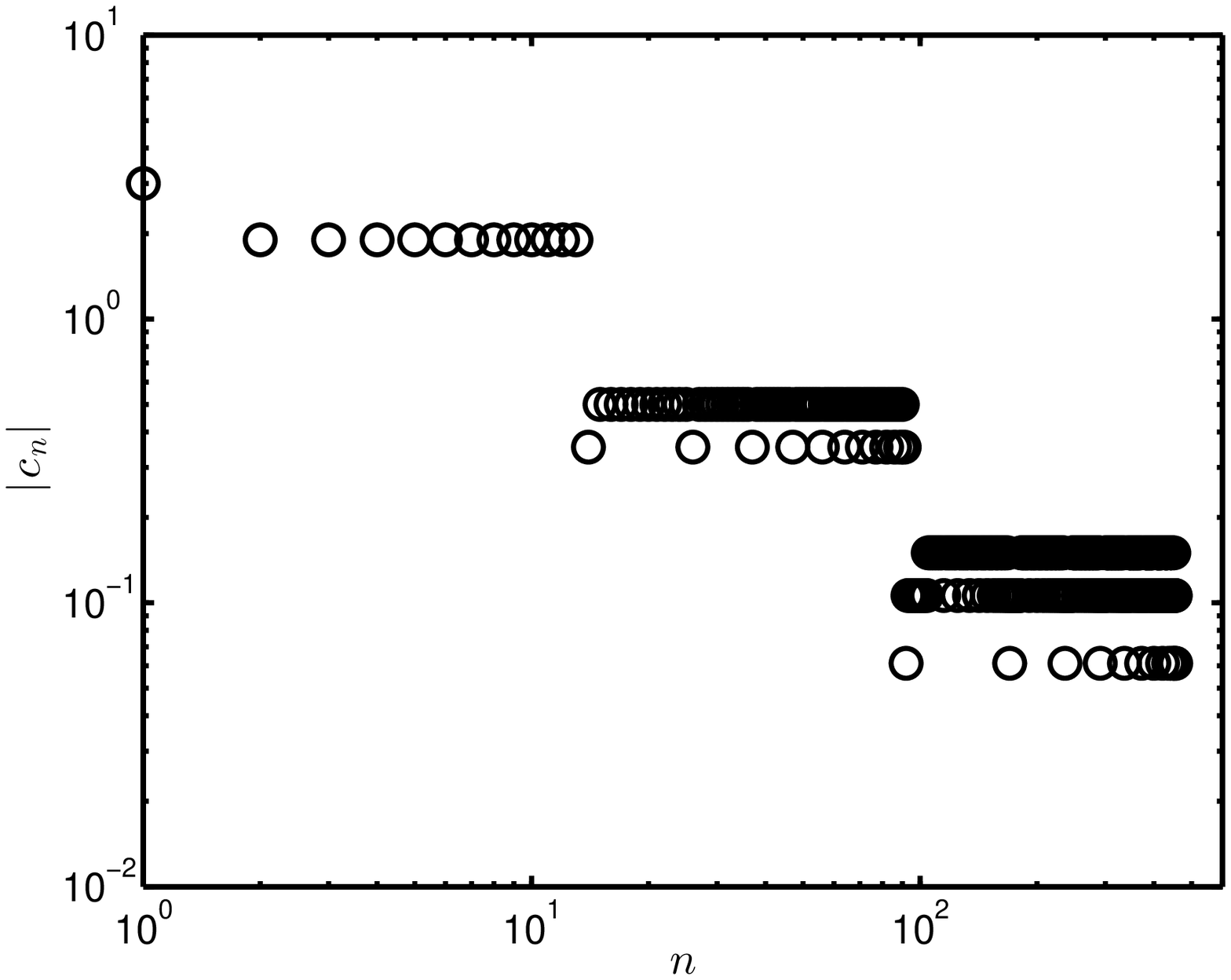}
\includegraphics[width=3in]{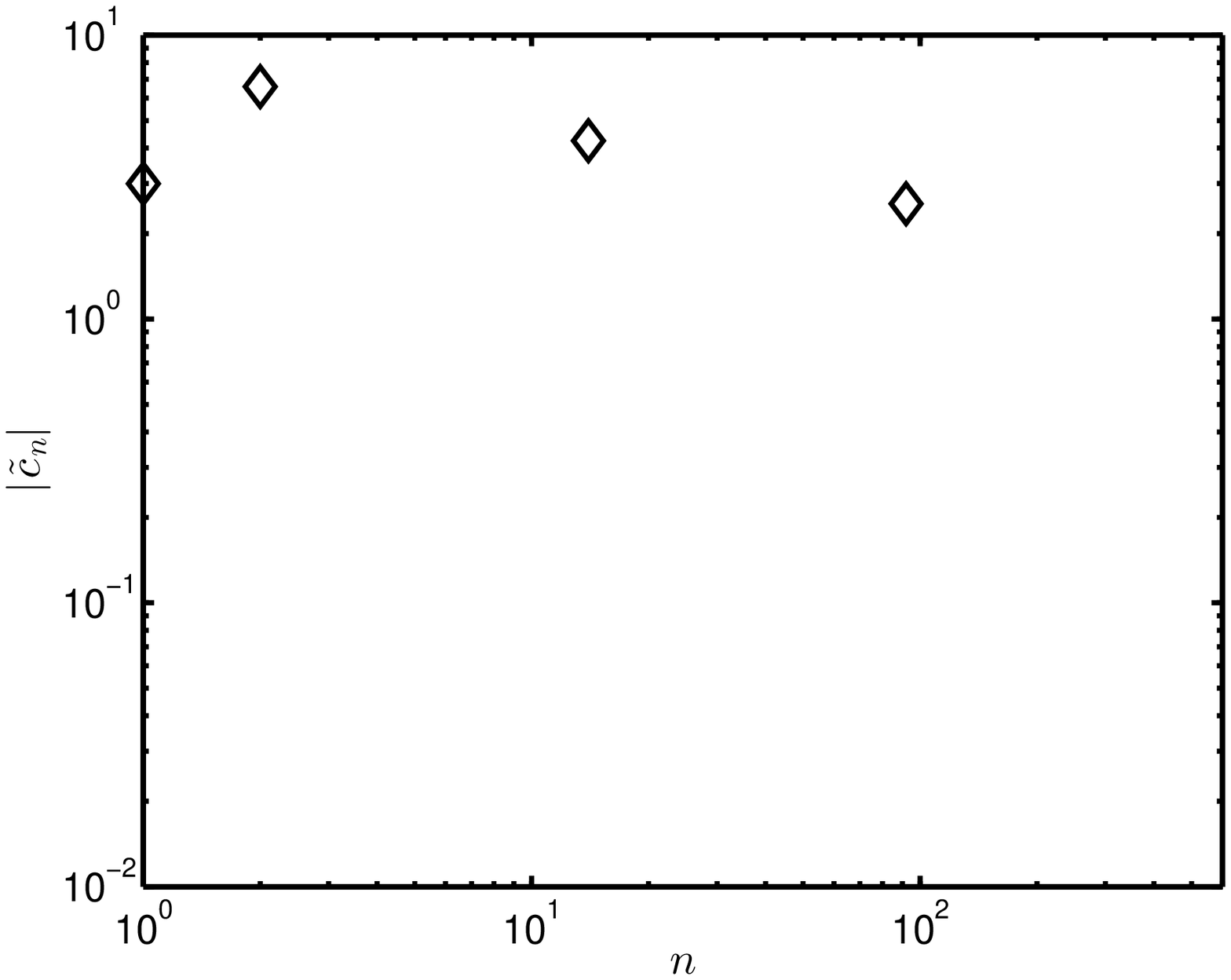}
\caption{Magnitude of the gPC coefficients for the example function with 
equally important random variables. (Left) Magnitude of $\bm c_n$. (Right) 
Magnitude of $\tilde{\bm c}_n$ of a randomly chosen replicate computed by
rotated $\ell_1$ with $9$ iterations and $M=180$ ($M/N\approx 0.4$). These 
calculations were performed with dimension $d=12$ and the number of unknowns
$N=455$.}
\label{fig:ex1_coef}
\end{figure}

In this test, we set $d=12$ (hence, $N=455$ for $P=3$) and demonstrate the 
effectiveness of our new method. The integrals for calculating the $L_2$ error 
are computed by a level $4$ sparse grid method, hence they are exact. The 
relative error of $\ell_1$ minimization and OMP are presented in 
Fig.~\ref{fig:ex1_l1}. Clearly, the standard $\ell_1$ minimization and OMP are
not effective as the relative error is close to $100\%$ even when $M/N > 0.4$.
Also, the re-weighted $\ell_1$ does not help in this case. However, our new 
iterative rotation demonstrates much better accuracy, especially when $M$ is
large. As demonstrated in Fig.~\ref{fig:ex1_coef} the iterative rotation creates
a much sparser representation of $u$, hence the efficiency of compressive 
sensing method is substantially enhanced. We notice that the accuracy increases
as more iterations are included. However, the improvement from $6$ iterations to
$9$ iterations is less significant as that from $3$ iterations to $6$ 
iterations, especially for the OMP-based iterative rotation method. In general,
the improvement afforded by iterative rotation becomes small after $3$ 
iterations.

This contrived example demonstrates that our new method is capable of enhancing
the sparsity of the Hermite polynomial expansion, even with a very inaccurate
$\tilde{\bm c}^{(0)}$ in Step 2 of Algorithm \ref{algo:cs_rot} when other 
methods fail.

\subsection{Example function with high compressibility} 
\label{subsec:ex2}
\begin{figure}[t]
\centering
\includegraphics[width=3in]{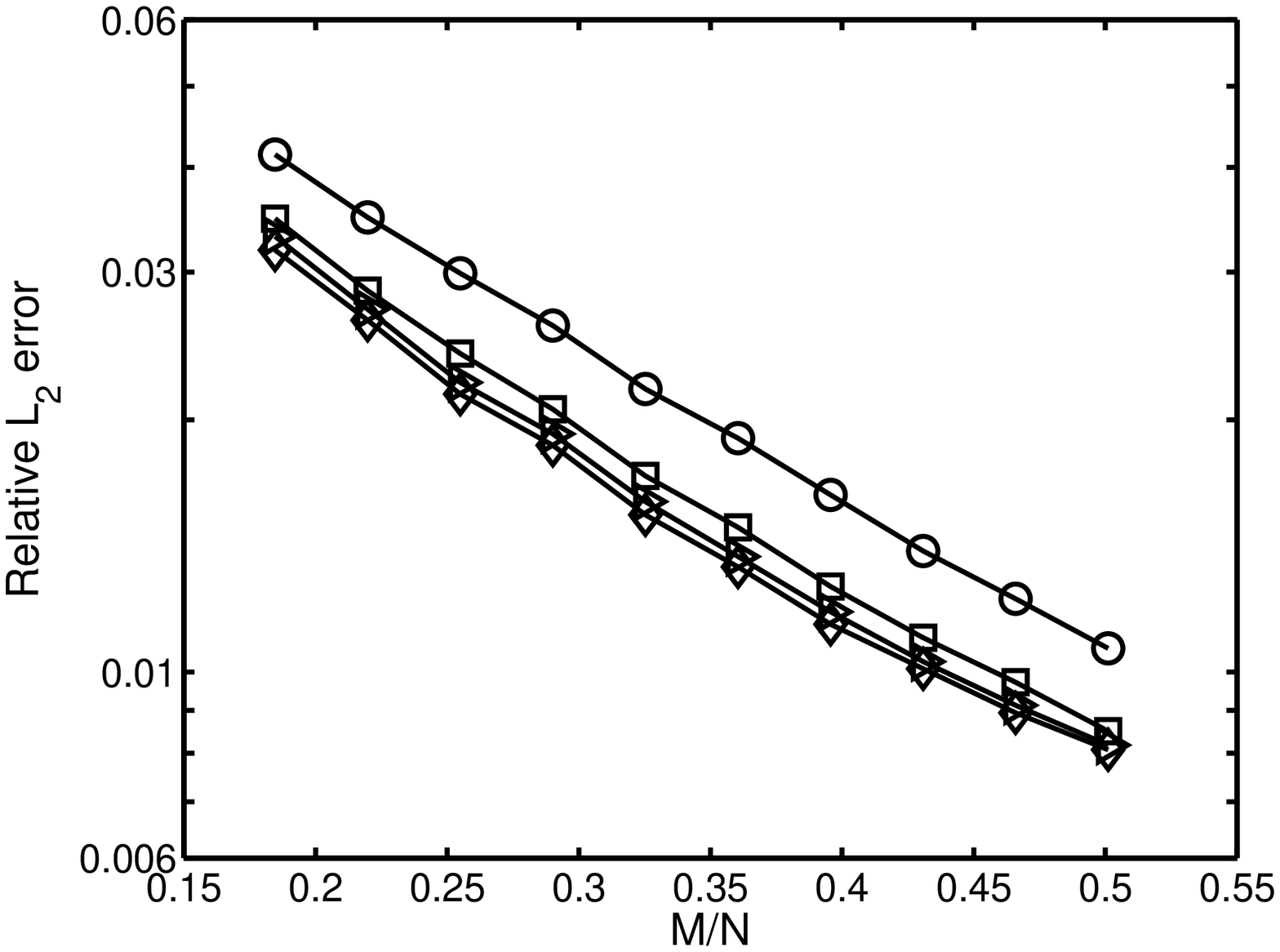}
\includegraphics[width=3in]{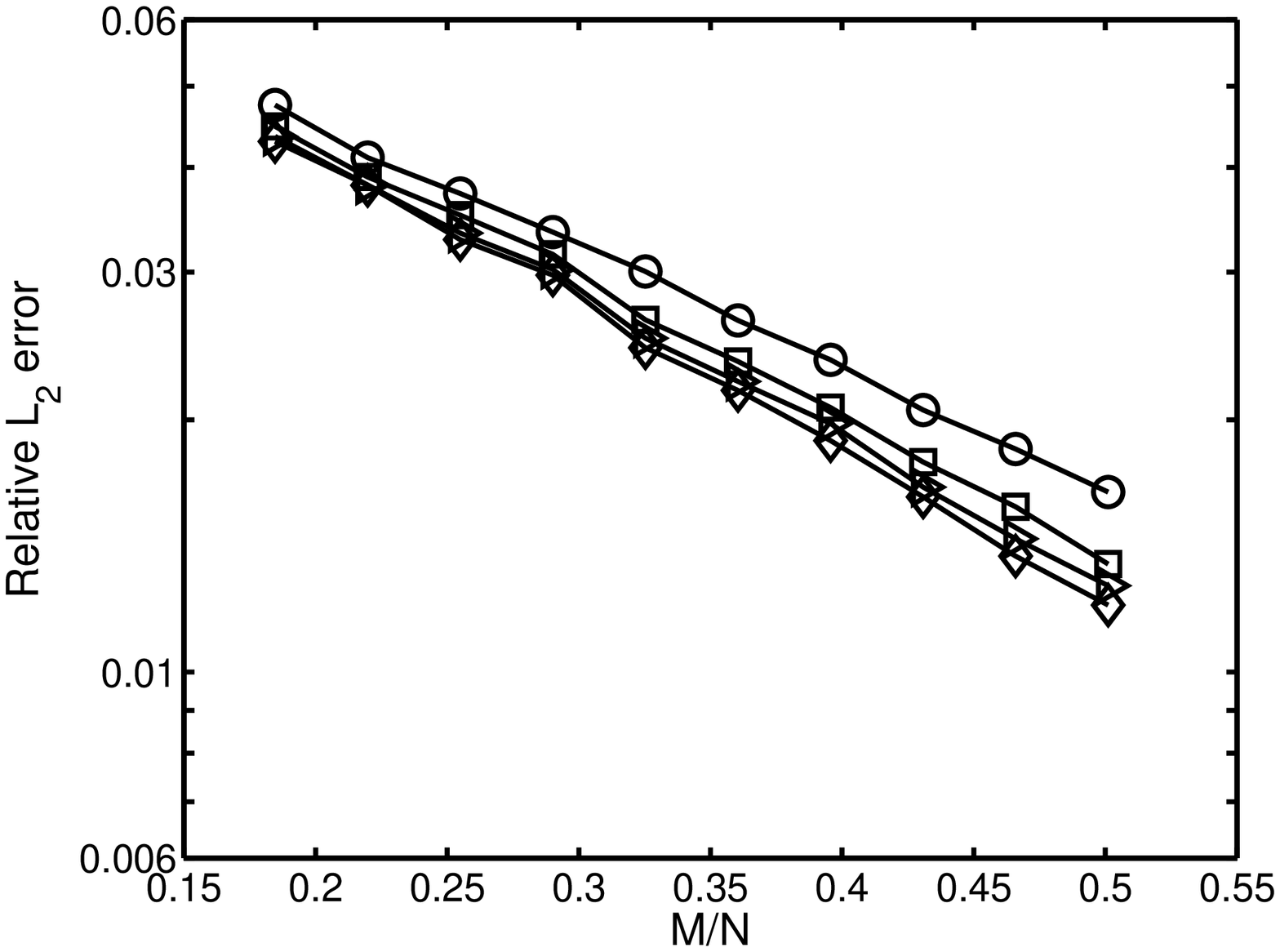}
\caption{Results for the example function with high compressibility. Relative 
$L_2$ error for different numbers of iterations for $\ell_1$ minimization (left)
and the OMP method (right). ``$\circ$": standard $\ell_1$ (left) or standard OMP
(right), ``$\Box$": $1$ iteration, ``$\triangleright$": $2$ iterations, 
``$\diamond$": $3$ iterations. These calculations were performed with dimension
$d=12$ and number of unknowns $N=455$.}
\label{fig:ex2_iter}
\end{figure}

Consider the following function:
\begin{equation}\label{eq:ex2}
u(\bx) = \sum_{|\ba|=0}^P c_{\ba}\psi_{\ba}(\bx)
       = \sum_{n=1}^Nc_n\psi_n(\bx), \quad
\bx =  (\xi_1,\xi_2,\cdots,\xi_{d}),
\end{equation}
where, $\psi_{\ba}$ are normalized multivariate Hermite polynomials, 
$d=12, P=3, N=455$, and the coefficients $c_n$ are chose as uniformly 
distributed random numbers,
\begin{equation}
c_n = \zeta/n^{1.5}, \quad \zeta\sim \mathcal{U}[-1,1].
\end{equation}
For this example, we generate $N$ samples of $\zeta$: 
$\zeta^1,\zeta^2,\cdots,\zeta^N$ then divide them by $n^{1.5}, n=1,2,\cdots,N$
to obtain a random ``compressible signal" $\bm c$. The integrals for the 
relative error are computed by a level-$4$ sparse grid method and are therefore
exact. Figure \ref{fig:ex2_iter} shows the the relative error with different 
numbers of iterations ($1$-$3$) for the $\ell_1$ minimization and OMP methods.
Our new iterative rotation method improves the accuracy of the gPC approximation 
for both methods. As before, benefit of increased iterations drops sharply near
$L=3$. Therefore, in the remainder of this paper we use $L=3$ iterations unless 
otherwise noted.

\begin{figure}[t]
\centering
\includegraphics[width=3in]{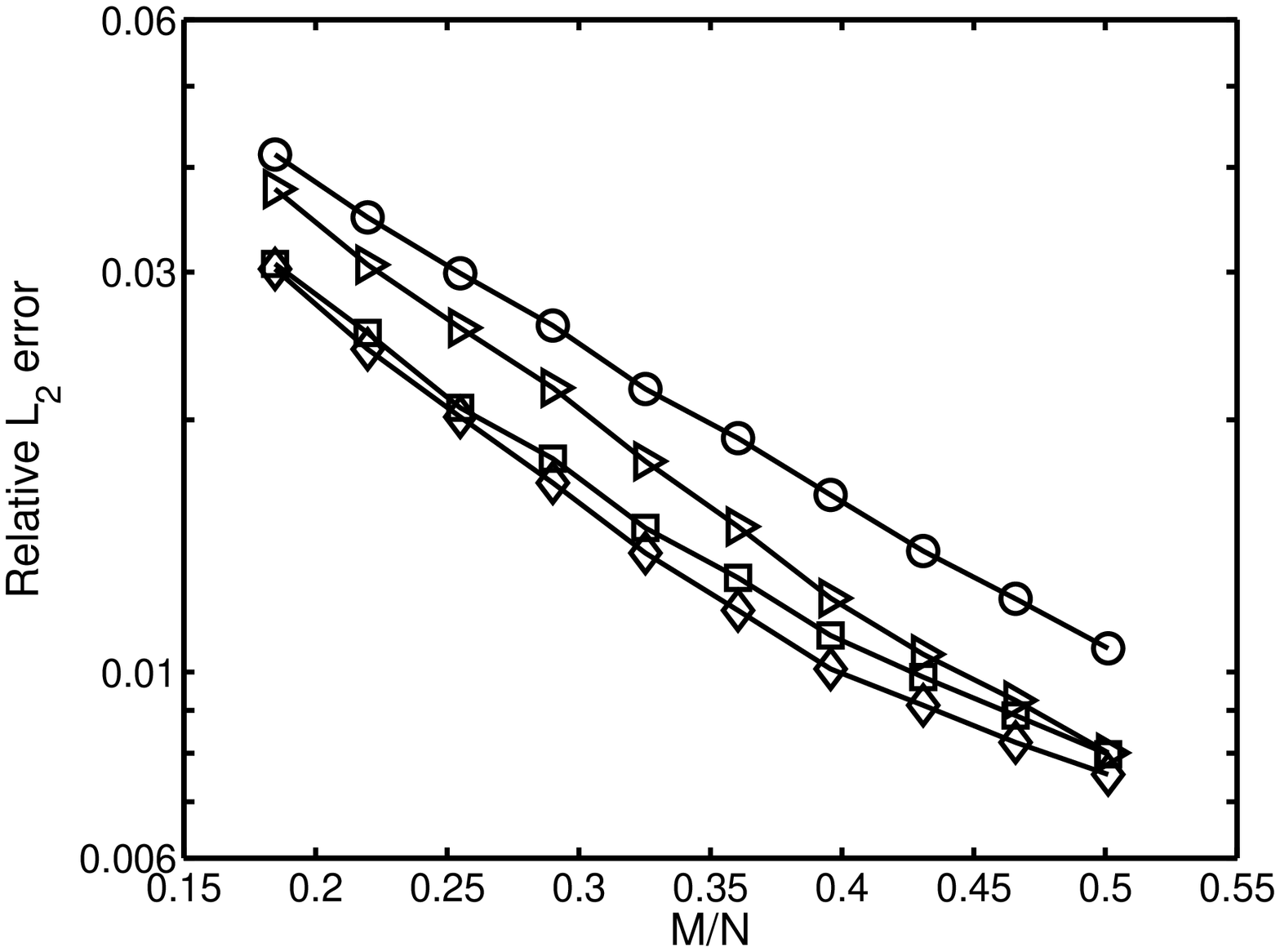}
\includegraphics[width=3in]{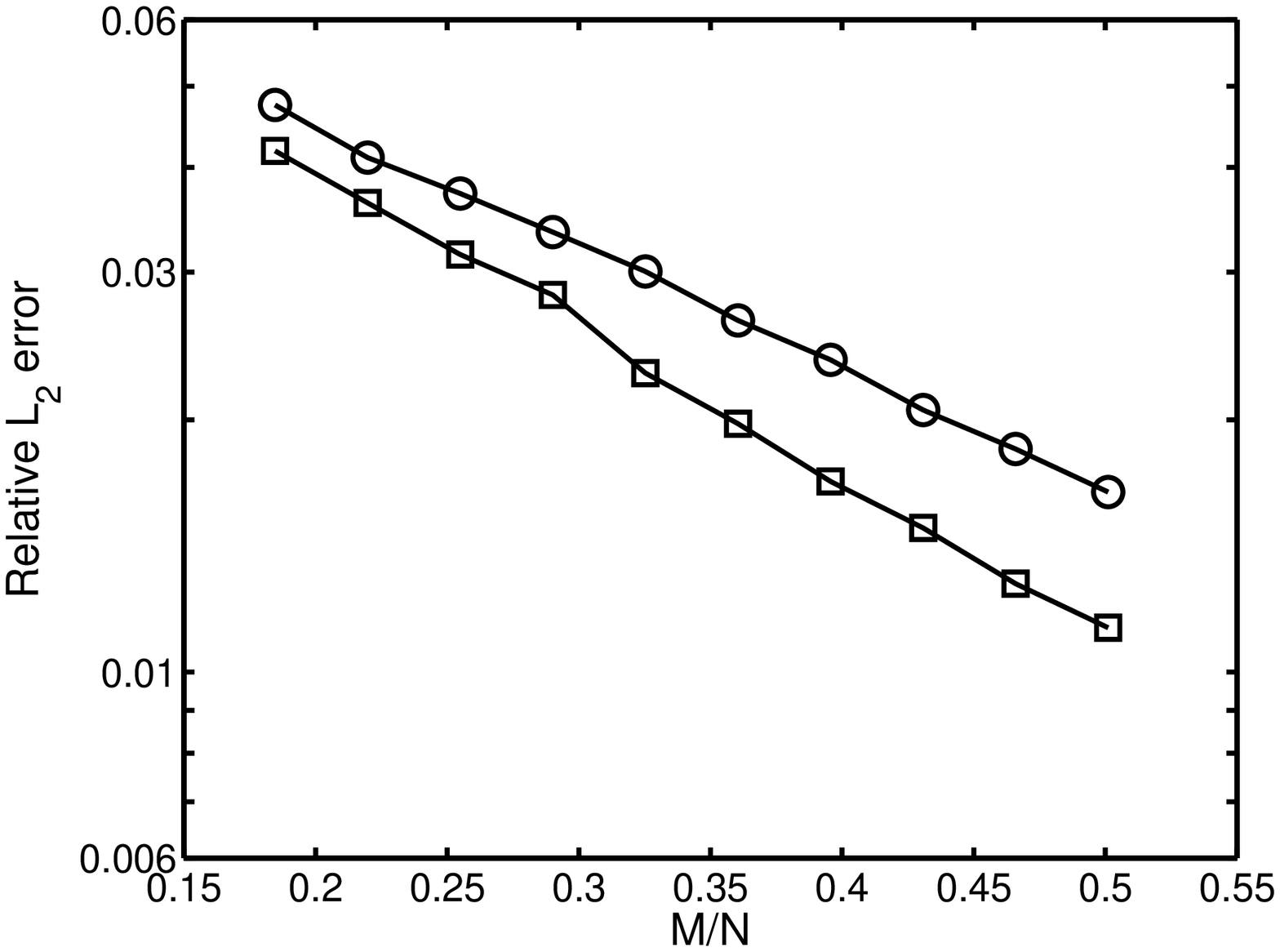}
\caption{Results for the example function with high compressibility. (Left)
Comparison with $\ell_1$ methods. ``$\circ$": standard $\ell_1$,
``$\triangleright$": re-weighted $\ell_1$, ``$\Box$": rotated $\ell_1$,
``$\diamond$": re-weighted and iteratively rotated $\ell_1$. (Right) Comparison
with OMP methods. ``$\circ$": OMP, ``$\Box$": rotated OMP. These calculations
were performed with dimension $d=12$ and number of unknowns $N=455$.}
\label{fig:ex2_l1}\label{fig:ex2_omp}
\end{figure}

\begin{figure}[t]
\centering
\includegraphics[width=3in]{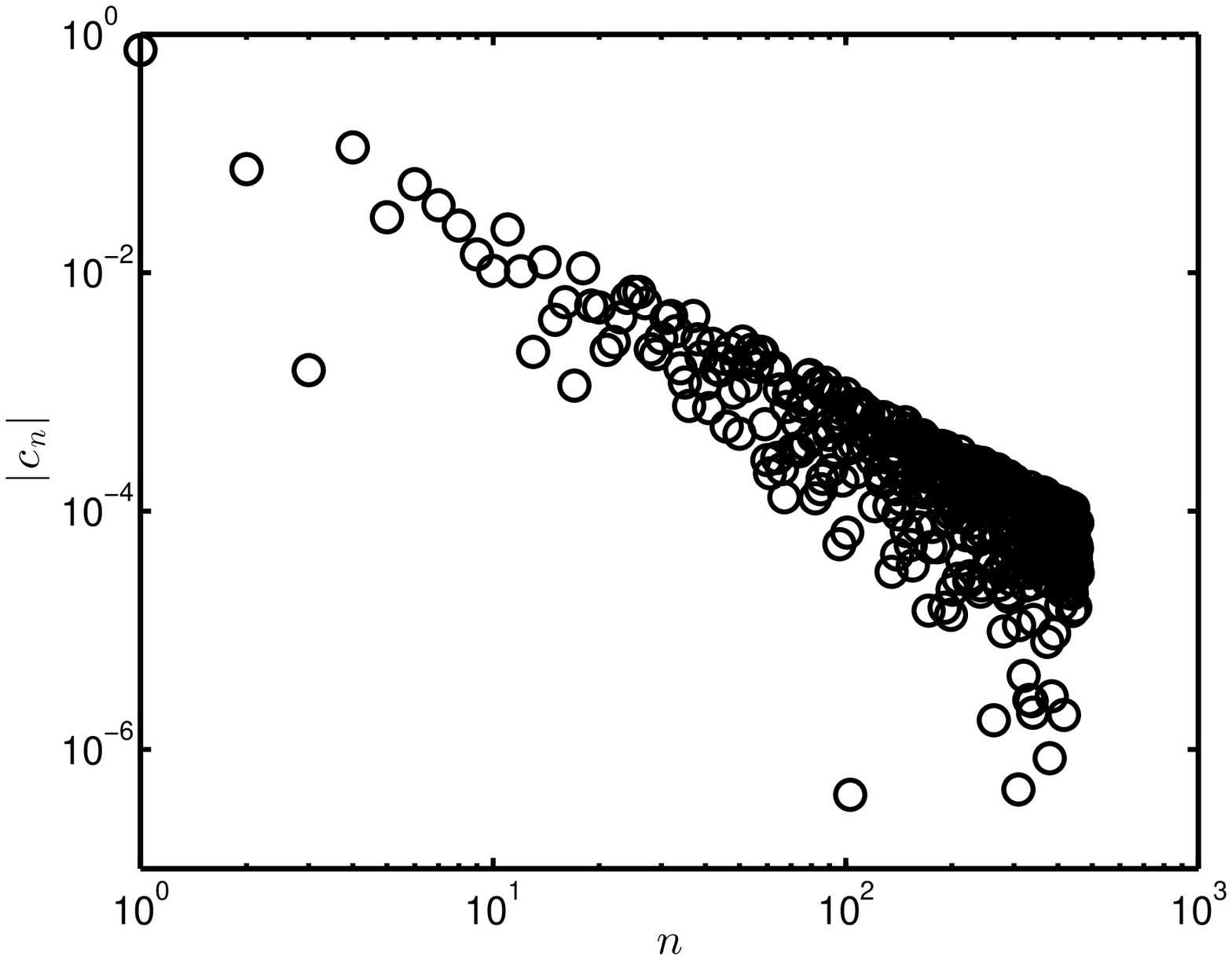}
\includegraphics[width=3in]{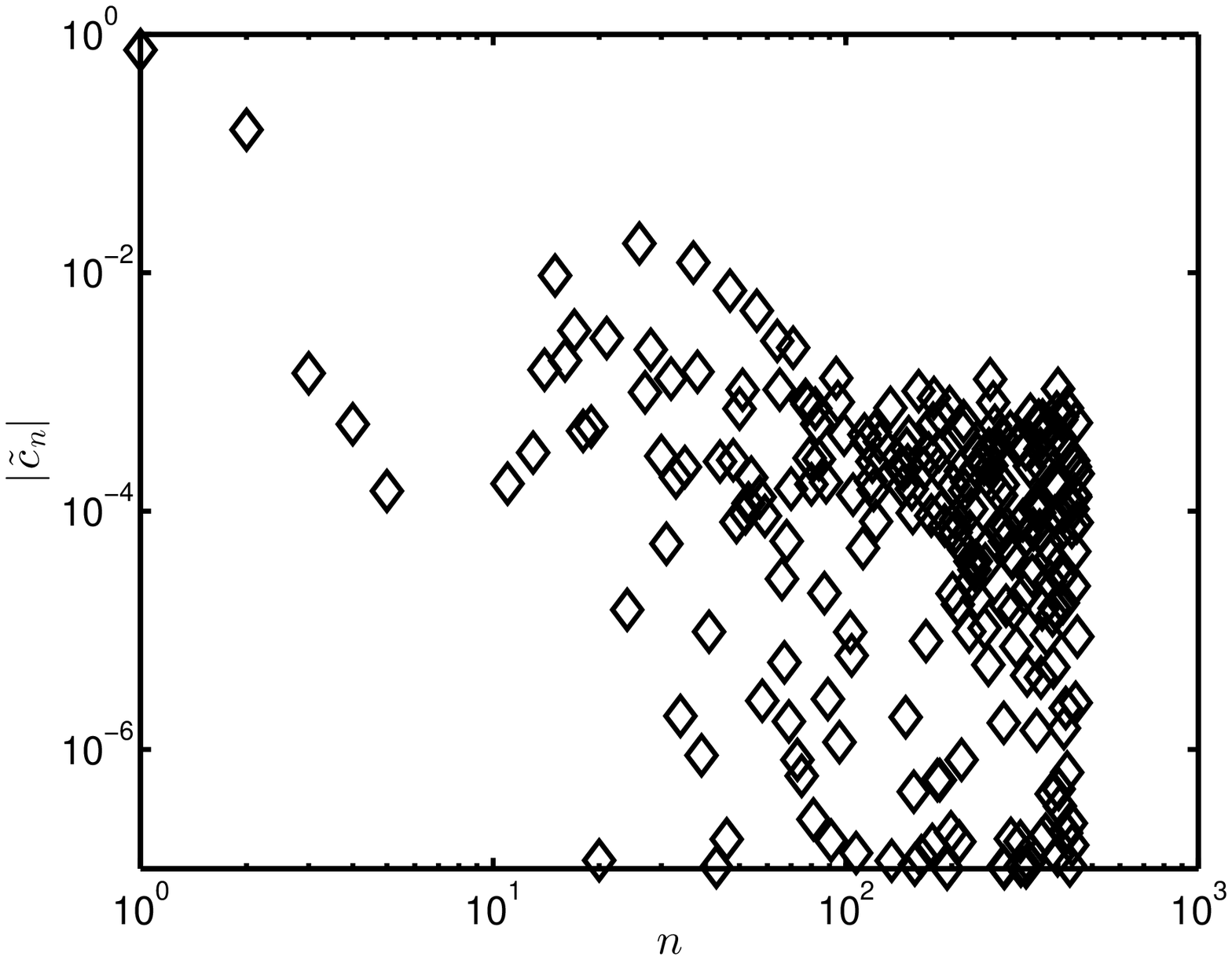}
\caption{Magnitude of the gPC coefficients for example function with high 
compressibility. (Left) Magnitude of $\bm c_n$. (Right) Magnitude of 
$\tilde{\bm c}_n$ of a randomly chosen replicate computed by re-weighted and 
iteratively rotated $\ell_1$ with $M=180$ ($M/N\approx 0.4$). These 
calculations were performed with dimension $d=12$ and the number of unknowns
$N=455$.}
\label{fig:ex2_coef}
\end{figure}

Figure \ref{fig:ex2_l1} shows results obtained by applying our iterative 
rotation technique to the re-weighted $\ell_1$ approach using $L=3$ iterations.
The results for the iterative rotation approach for OMP are also presented in
Figure \ref{fig:ex2_l1}. For all methods, introduction of the iterative rotation
approach improves the results. A comparison of the sparsity of $\bm c$ and
$\tilde{\bm c}$ is presented in Fig.~\eqref{fig:ex2_coef}. The main improvement
is that the number of coefficients with magnitude larger than $0.01$ is 
decreased. Also, $|c_n|$ cluster around the line $|c_n|=1/n^{1.5}$ as we set 
them in this way, while many $|\tilde c_n|$ are much below this line especially 
when $n$ is large.

\subsection{Example elliptic differential equation}
\label{subsec:ex3}
\begin{figure}[t]
\centering
\includegraphics[width=3in]{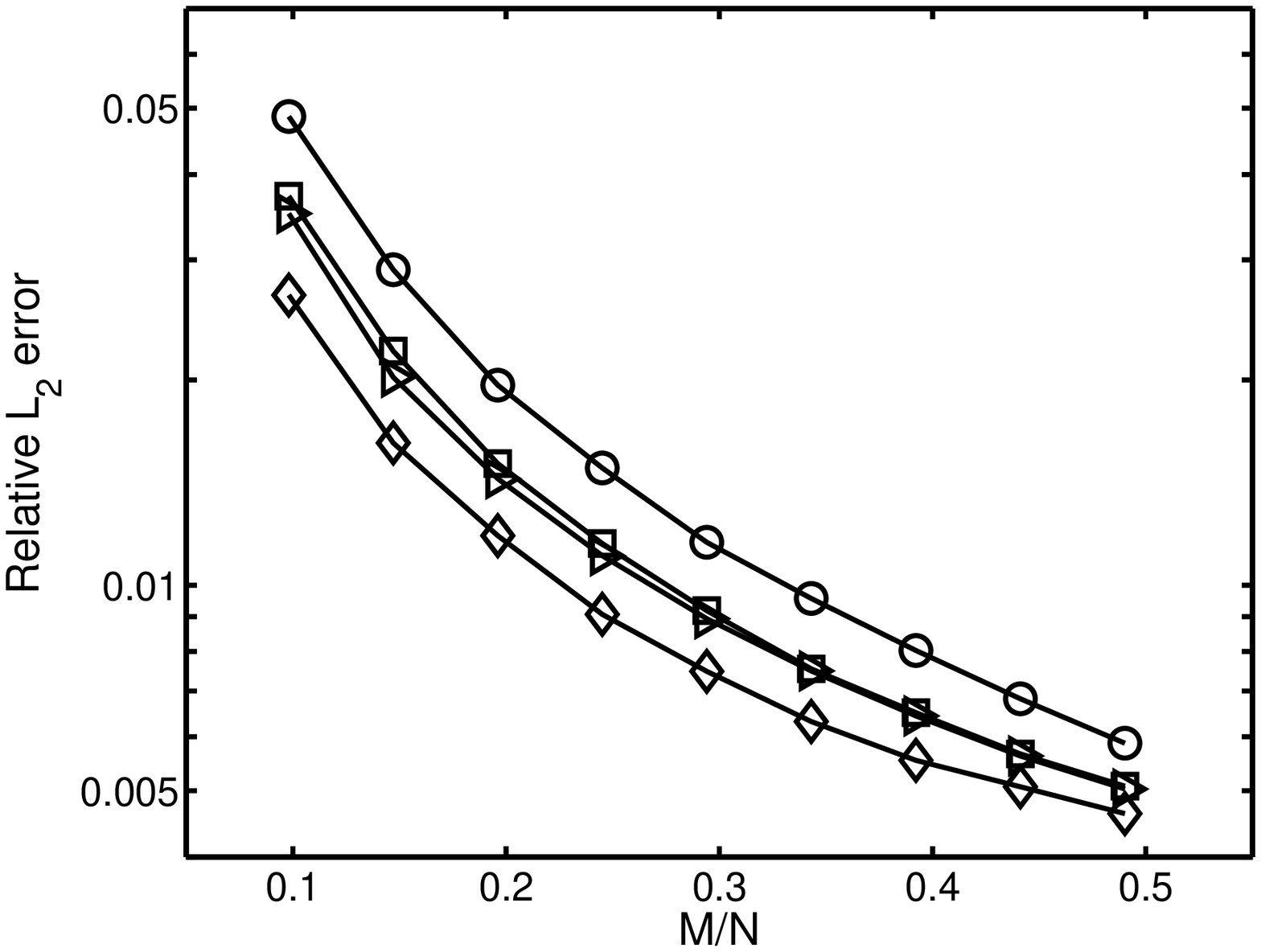}
\includegraphics[width=3in]{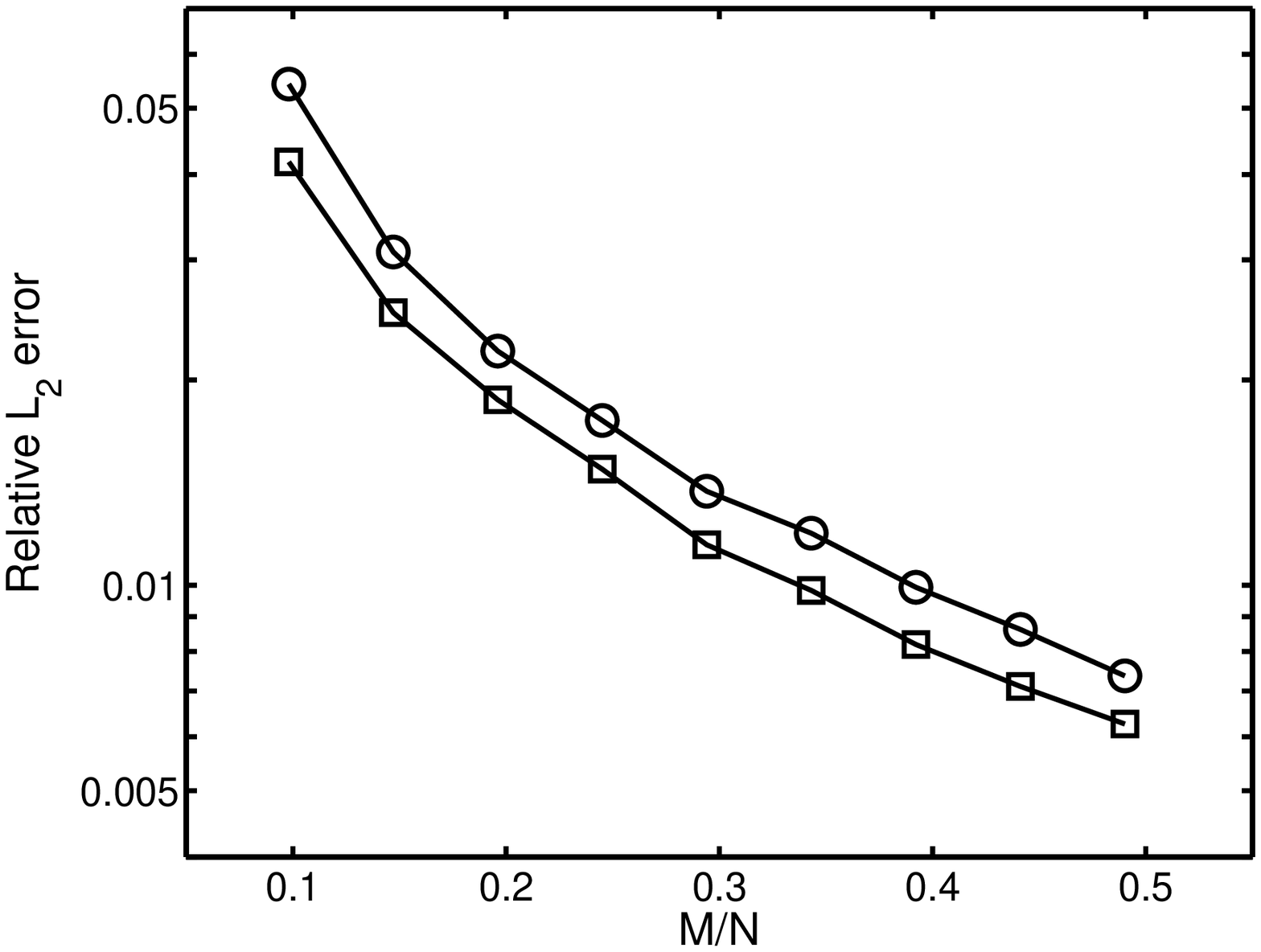}
\caption{Results for the example elliptic differential equation. 
(Left) Comparison with $\ell_1$ methods.  ``$\circ$": standard $\ell_1$,
``$\triangleright$": re-weighted $\ell_1$, ``$\Box$": rotated $\ell_1$, 
``$\diamond$": re-weighted and iteratively rotated $\ell_1$. (Right) Comparison
with OMP methods. ``$\circ$": standard OMP, ``$\Box$": rotated OMP. These 
calculations were performed with dimension $d=15$ and the number of unknowns
$N=816$.}
\label{fig:ex3_l1} \label{fig:ex3_omp}
\end{figure}

\begin{figure}[t]
\centering
\includegraphics[width=3in]{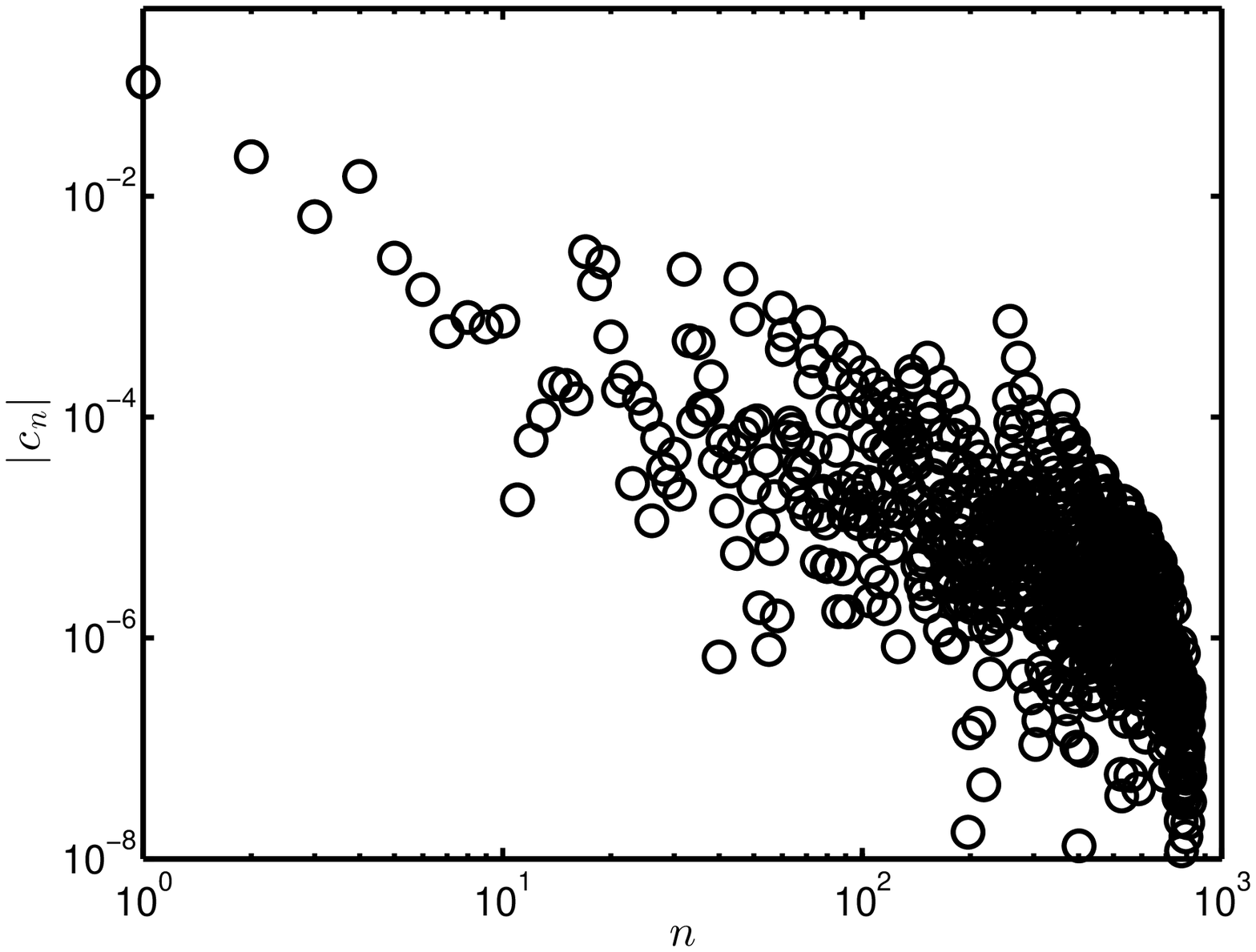}
\includegraphics[width=3in]{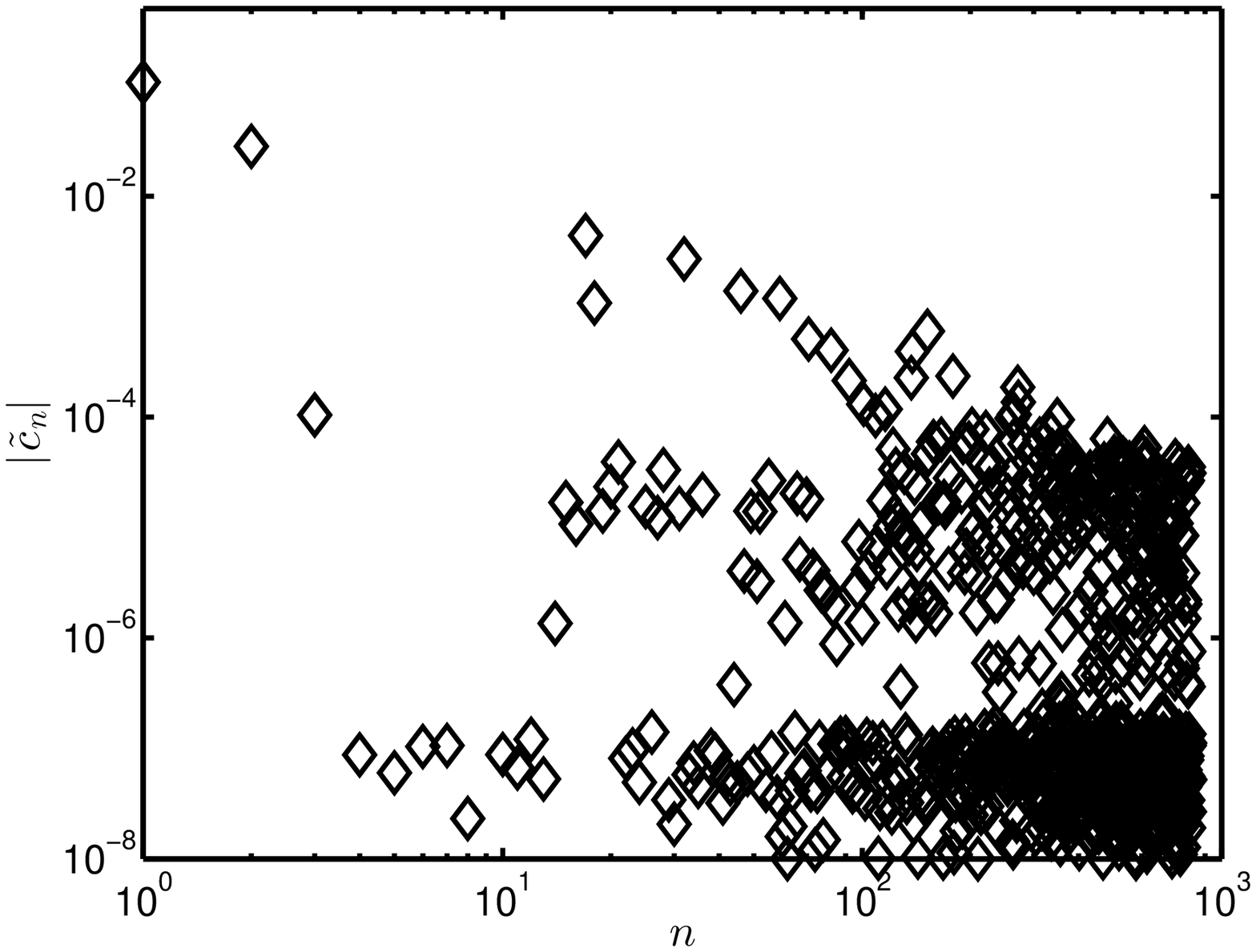}
\caption{Magnitude of the gPC coefficients for example elliptic differential
equation. (Left) Magnitude of $\bm c_n$.
(Right) Magnitude of $\tilde{\bm c}_n$ of a randomly chosen replicate computed
by re-weighted and iteratively rotated $\ell_1$ with $M=240$ ($M/N\approx 0.3$).
These calculations were performed with dimension $d=15$ and the number of 
unknowns $N=816$.}
\label{fig:ex3_coef}
\end{figure}

Next we consider a one-dimensional elliptic differential equation with a random 
high-order coefficient:
\begin{equation}\label{eq:ellip}
\begin{aligned}
-\frac{d}{dx} \left( a(x;\bx)\frac{d u(x;\bx)}{dx} \right) = 1, & \quad x \in (0,1) \\
u(0) = u(1) = 0, &
\end{aligned}
\end{equation}
where $a(x;\bx)$ is a log-normal random field based on Karhunen-Lo\`eve (KL) 
expansion:
\begin{equation}\label{eq:kl}
a(x;\bx) = a_0(x) +
\exp\left(\sigma\sum_{i=1}^d\sqrt{\lambda_i}\phi_i(x)\xi_i\right), 
\end{equation}
where $\{\xi_i\}$ are i.i.d.~standard Gaussian random variables,
$\{\lambda_i\}_{i=1}^d$, and $\{\phi_i(x)\}_{i=1}^d$ are the largest eigenvalues
and corresponding eigenfunctions of the exponential covariance kernel:
\begin{equation}\label{eq:exp_kernel}
C(x,x') = \exp\left(\dfrac{|x-x'|}{l_c}\right).
\end{equation}
In the KL expansion, $\lambda_i$ denotes the eigenvalue of the covariance kernel
$C(x,x')$ instead of entries of $\Lambda$ in Eq.~\eqref{eq:grad}. The value of 
$\lambda_i$ and the analytical expressions for $\phi_i$ are available in the 
literature \cite{JardakSK02}. In this example, we set 
$a_0(x)\equiv 0.1, \sigma=0.5, l_c=0.2$ and $d=15$. With this setting, 
$\sum_{i=1}^d\lambda_i > 0.93\sum_{i=1}^{\infty}\lambda_i$. For each input 
sample $\bx^q$, $a$ and $u$ only depend on $x$ and the solution of the 
deterministic elliptic equation can be obtained as \cite{YangK13}:
\begin{equation}\label{eq:ellip_sol}
u(x) = u(0) + \int_0^x \dfrac{a(0)u(0)'-y}{a(y)}\dif y.
\end{equation}
By imposing the boundary condition $u(0)=u(1)=0$, we can compute $a(0)u(0)'$ as
\begin{equation}\label{eq:ellip_const}
a(0)u(0)' = \left(\int_0^1 \dfrac{y}{a(y)}\dif y\right) \Big /
            \left(\int_0^1 \dfrac{1}{a(y)}\dif y\right).
\end{equation}
The integrals in Eqs.~\eqref{eq:ellip_const} and \eqref{eq:ellip_sol} must be
obtained by highly accurate numerical integration. For this example, we choose 
the quantity of interest to be $u(x;\bx)$ at $x=0.35$. We aim to build a 
3rd-order Hermite polynomial expansion which includes $N=816$ basis functions. 
The relative error is approximated by a level-$6$ sparse grid method. 
Figure \ref{fig:ex3_l1} shows that accuracy of the re-weighted $\ell_1$ ($3$
iterations) and the iteratively rotated $\ell_1$ ($L = 3$ iterations) method 
are very close in this case. Figure \ref{fig:ex3_omp} shows the results of the
iterative rotation process applied to the OMP method. In all cases, the 
incorporation of iterative rotation improves the performance of the other 
methods. A comparison of $\bm c$ and $\tilde{\bm c}$ are presented in
Fig.~\ref{fig:ex3_coef}, which shows the improvement of the sparsity in the
similar manner as in example function with high compressibility in
Sec.~\ref{subsec:ex2}.

\subsection{Example Korteweg-de Vries equation} \label{subsec:ex4}
\begin{figure}[t]
\centering
\includegraphics[width=3in]{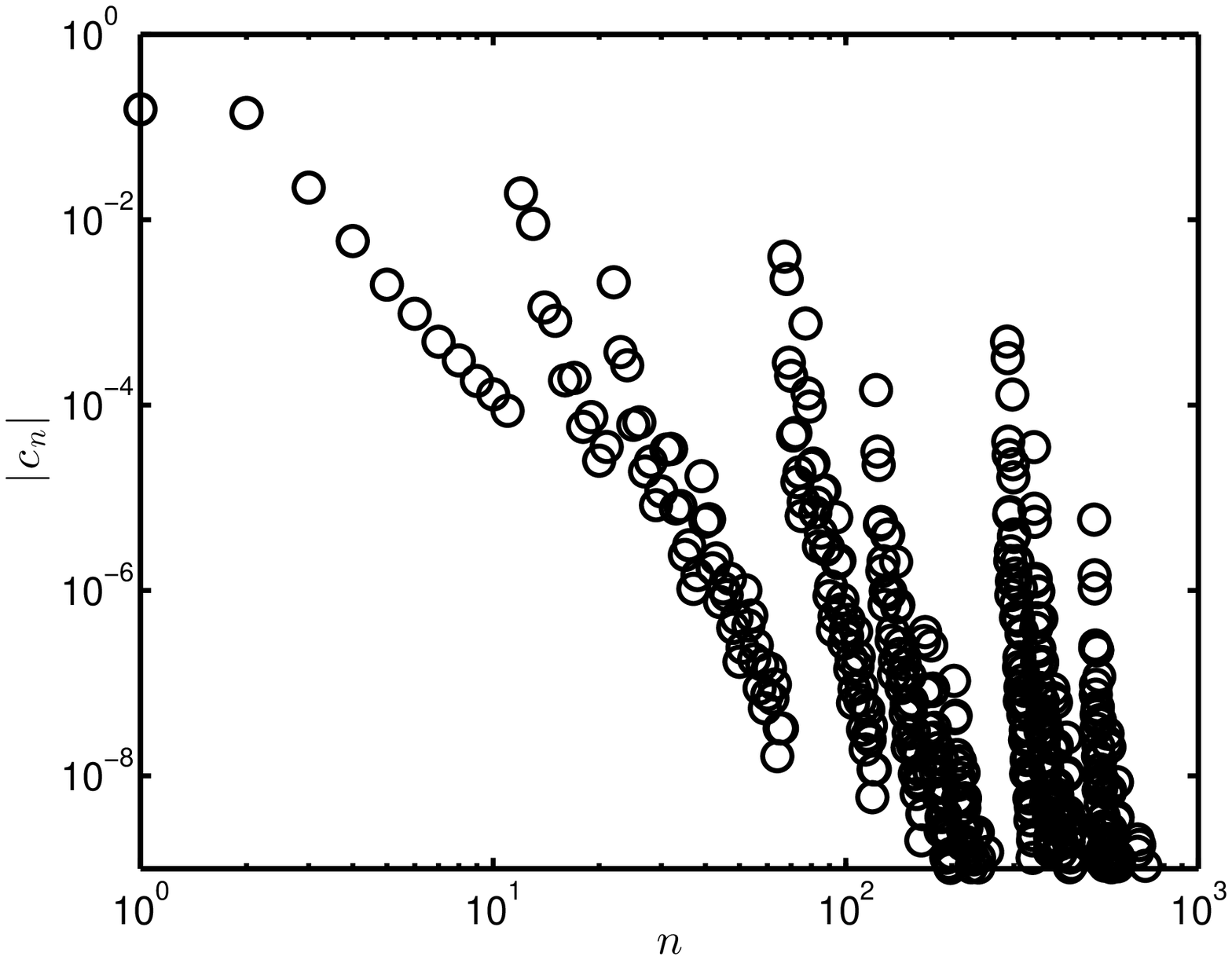}
\includegraphics[width=3in]{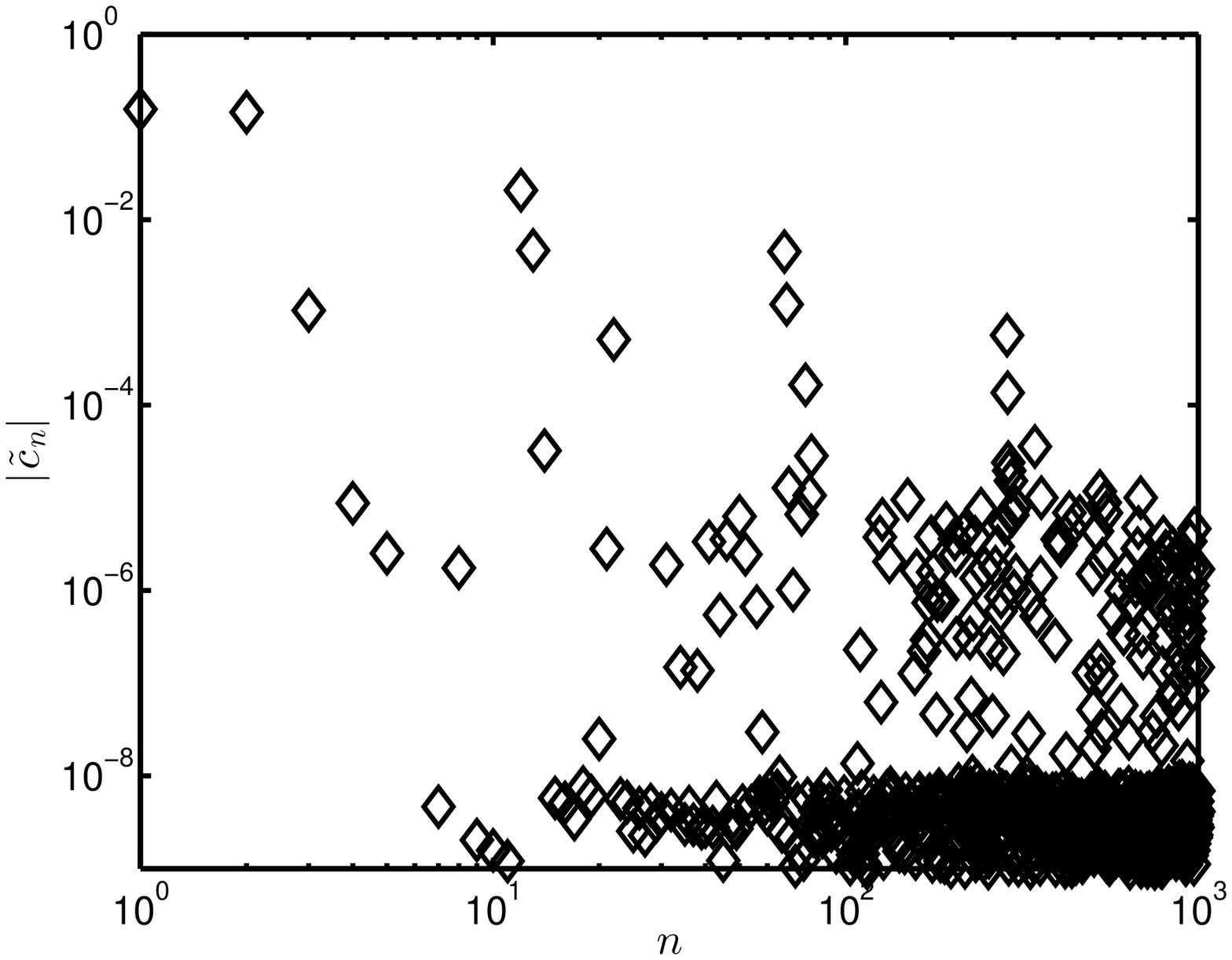}
\caption{Magnitude of the gPC coefficients for example KdV equation. (Left) 
Magnitude of $\bm c_n$. (Right) Magnitude of $\tilde{\bm c}_n$ of a randomly 
chosen replicate computed by re-weighted and iteratively rotated $\ell_1$ with 
$M=120$ ($M/N\approx 0.12$). These calculations were performed with dimension 
$d=10$ and the number of unknowns $N=1001$.}
\label{fig:ex4_coef}
\end{figure}

As an example application of our new method to a more complicated and nonlinear
differential equation, we consider the Korteweg-de Vries (KdV) equation with 
time-dependent additive noise \cite{LinGK06}:
\begin{equation}\label{eq:kdv}
\begin{aligned}
& u_t(x,t;\bx)-6u(x,t;\bx)u_x(x,t;\bx)+u_{xxx}(x,t;\bx)=f(t;\bx), 
  \quad x\in (-\infty,\infty), \\
& u(x,0;\bx) = -2 \sech^2(x).
\end{aligned}
\end{equation}
Defining
\begin{equation}
W(t;\bx) = \int_0^t f(y;\bx) \dif y,
\end{equation}
the analytical solution of Eq.~\eqref{eq:kdv} is
\begin{equation}
u(x,t;\bx)=
W(t;\bx)-2\sech^2\left(x-4t+6\int_0^t W(z;\bx)\dif z\right).
\end{equation}
We model $f(t;\bx)$ as a Gaussian random field represented by the following
KL expansion:
\begin{equation}
f(t;\bx) = \sigma\sum_{i=1}^d\sqrt{\lambda_i}\phi_i(t)\xi_i, 
\end{equation}
where $\sigma$ is a constant and $\{\lambda_i,\phi_i(t)\}_{i=1}^d$ are
eigenpairs of the exponential covariance kernel as in Eqs.~\eqref{eq:kl} and
\eqref{eq:exp_kernel}, respectively. In this problem, we set $l_c=0.25$ and 
$d=10$ ($\sum_{i=1}^d\lambda_i > 0.96\sum_{i=1}^{\infty}\lambda_i$). In this 
case, the exact one-soliton solution is
\begin{equation}\label{eq:kdv_sol}
u(x,t;\bx)=\sigma\sum_{i=1}^d\sqrt{\lambda_i}\xi_i\int_0^t\phi_i(y)\dif y-
2\sech^2\left(x-4t+6\sigma\sum_{i=1}^d\sqrt{\lambda_i}\xi_i
\int_0^t\int_0^{z}\phi_i(y)\dif y\dif z\right).
\end{equation}
Since an analytical expression for $\phi_i$ is available, we can compute the 
integrals in Eq.~\eqref{eq:kdv_sol} with high accuracy. Denoting
\begin{equation}
A_i = \sqrt{\lambda_i}\int_0^t\phi_i(y)\dif y,\quad 
B_i = \sqrt{\lambda_i}\int_0^t\int_0^{z}\phi_i(y)\dif y\dif z, \quad
i=1,2,\cdots,d,
\end{equation}
the analytical solution is
\begin{equation}\label{eq:kdv_sol2}
u(x,t;\bx)\big |_{x=6,t=1}=\sigma\sum_{i=1}^dA_i\xi_i
-2\sech^2\left(2+6\sigma\sum_{i=1}^d B_i\xi_i\right).
\end{equation}
The quantity of interest is chosen to be $u(x,t;\bx)$ at $x=6,t=1$ with 
$\sigma=0.1$, $P=4$, and the number of gPC basis functions $N=1001$.
For this example, the combined iterative rotation and re-weighted $\ell_1$ 
method outperforms all other approaches. However, unlike previous examples the 
non-rotated re-weighted $\ell_1$ method works better than our iteratively 
rotated unweighted method. This difference likely arises because $\bm c$ is
sparser in this case than in others, which makes re-weighted $\ell_1$ method 
more efficient. The pattern of sparsity in this case is different than previous
examples, hence the efficiency of identifying a good rotation matrix $\tensor A$
is different. A comparison of $\bm c$ and $\tilde{\bm c}$ are presented in
Fig.~\ref{fig:ex4_coef}, which shows the improvement of the sparsity 
{ by the iterative rotation method.}

\begin{figure}[t]
\centering
\includegraphics[width=3in]{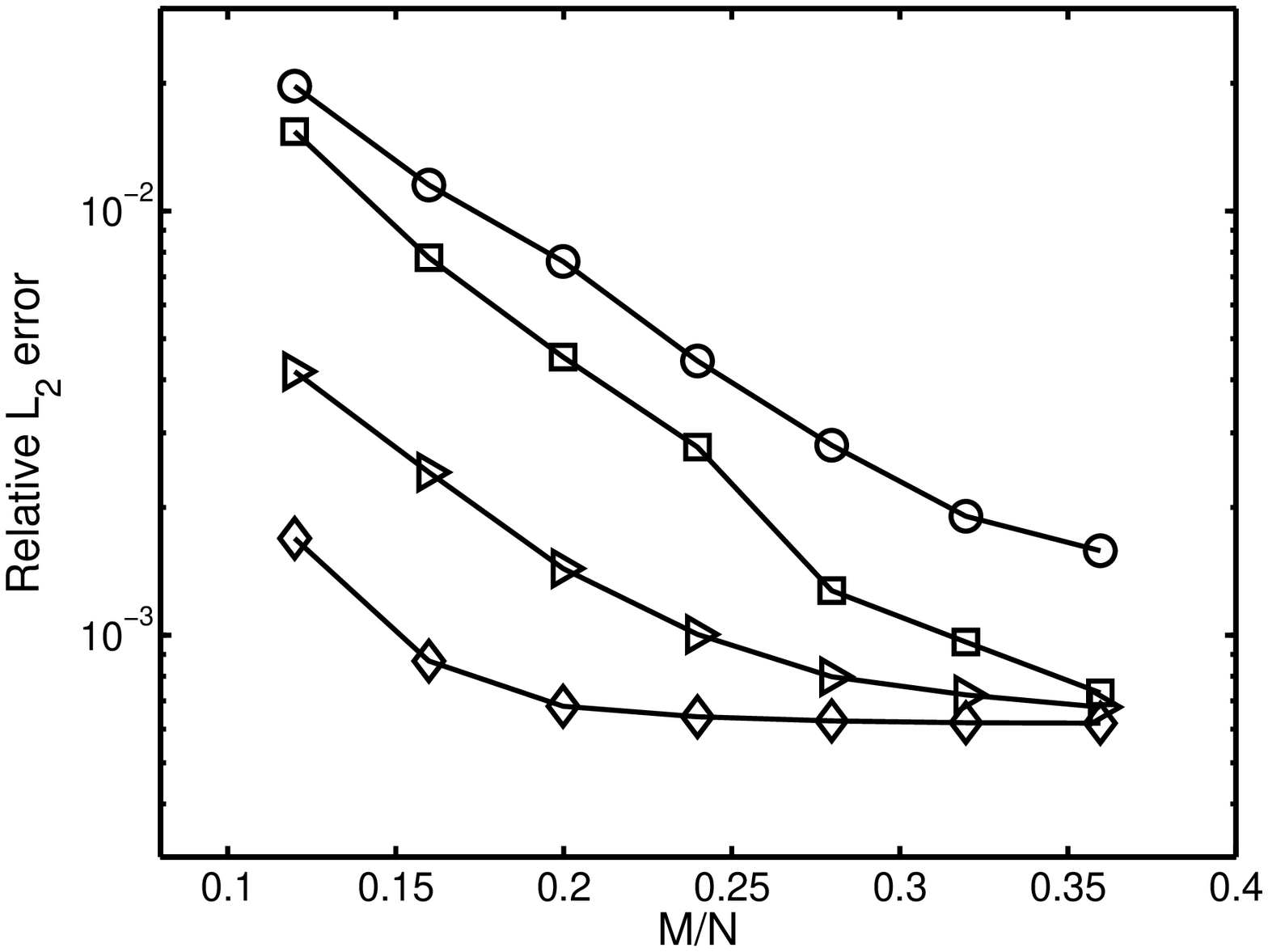}
\includegraphics[width=3in]{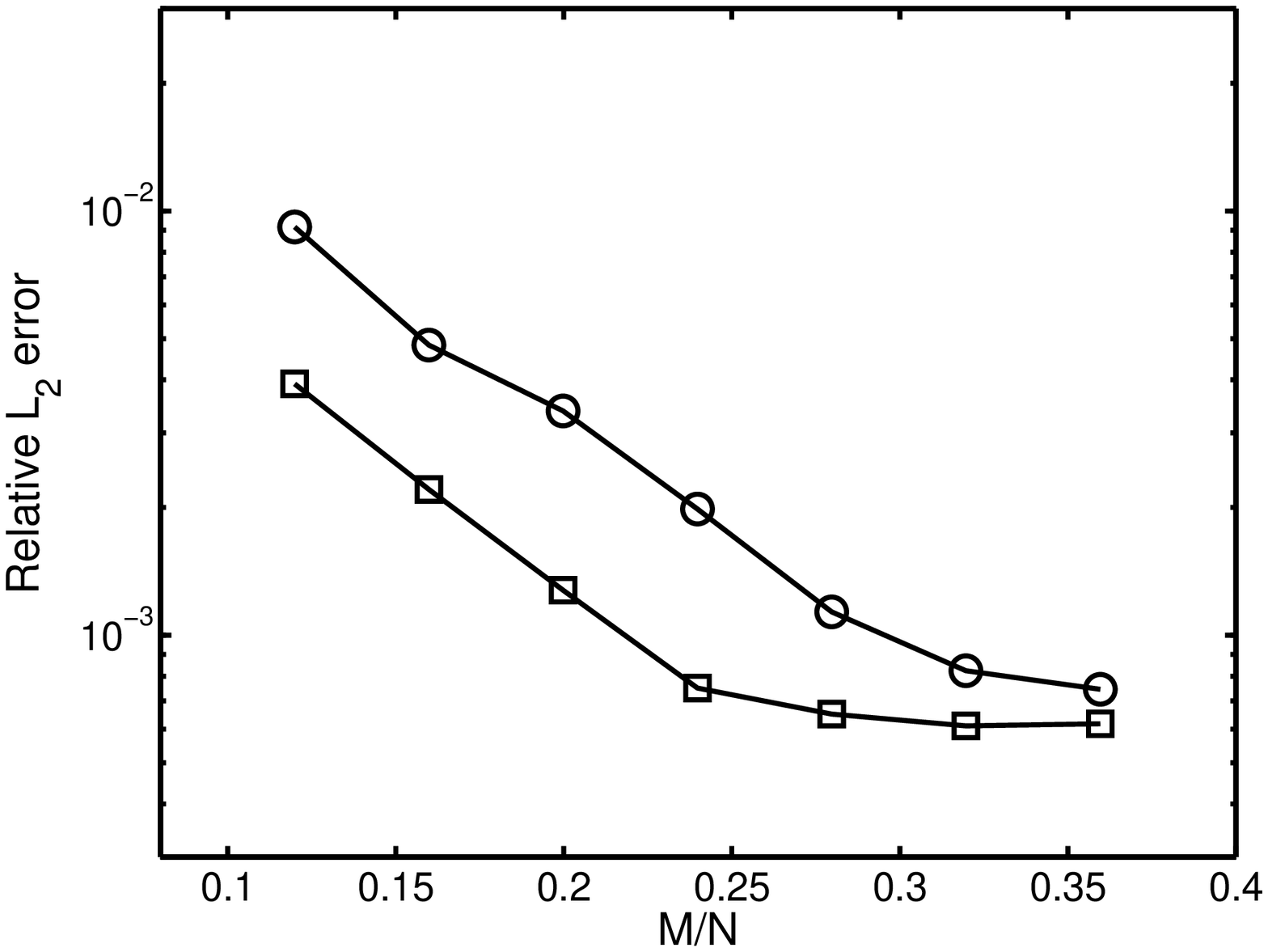}
\caption{Results for the example KdV equation. (Left) Comparison to $\ell_1$ 
methods.  ``$\circ$": standard $\ell_1$, ``$\triangleright$": re-weighted 
$\ell_1$, ``$\Box$": rotated $\ell_1$, ``$\diamond$": re-weighted and 
iteratively rotated $\ell_1$. (Right) Comparison to the OMP method.  $\circ$": 
standard OMP, ``$\Box$": rotated OMP. The dimension is $d=10$ and the number of
unknowns is $N=1001$.}
\label{fig:ex4_l1} \label{fig:ex4_omp}
\end{figure}

\subsection{Example high-dimensional function} \label{subsec:ex5}
\begin{figure}[t]
\centering
\includegraphics[width=3in]{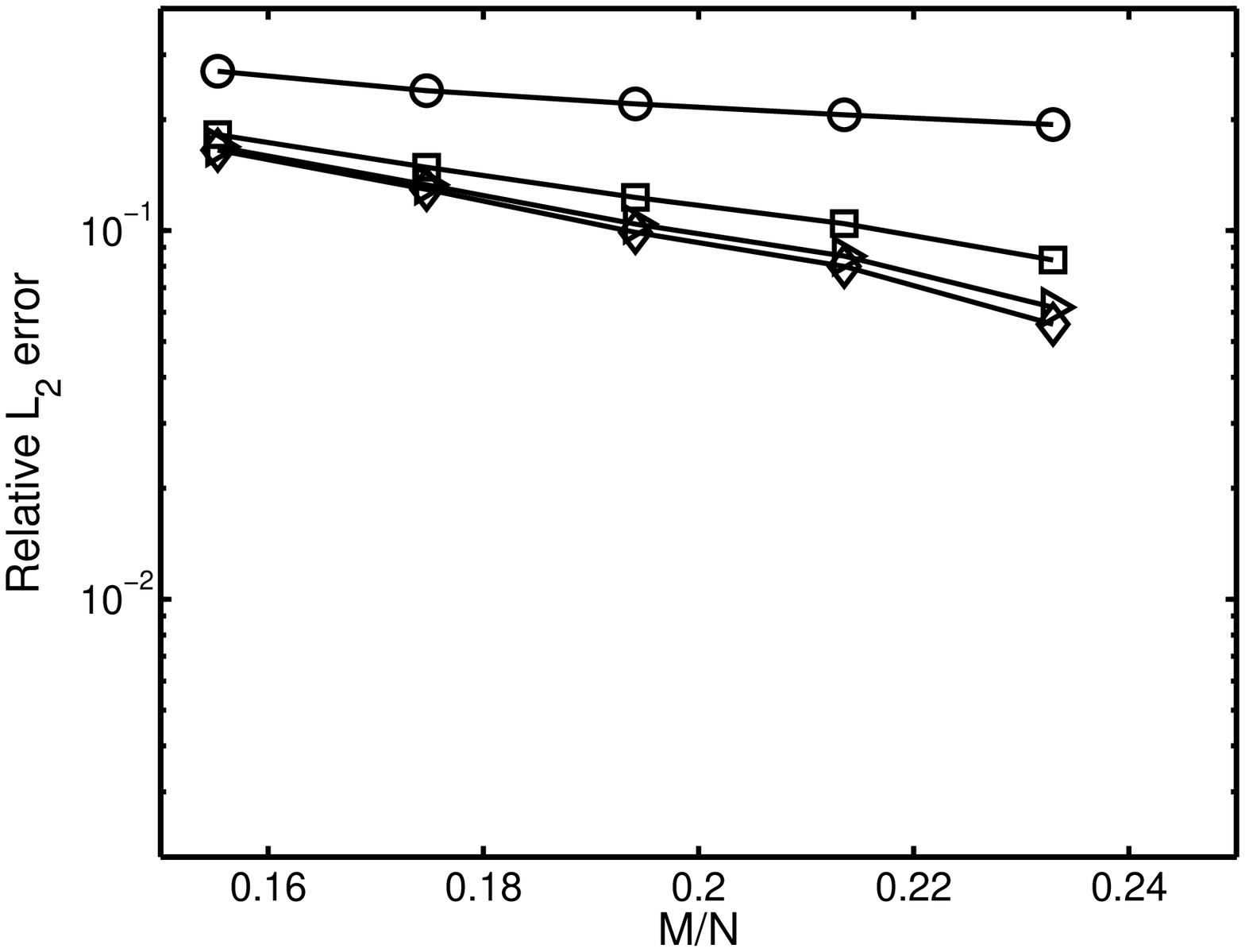}
\includegraphics[width=3in]{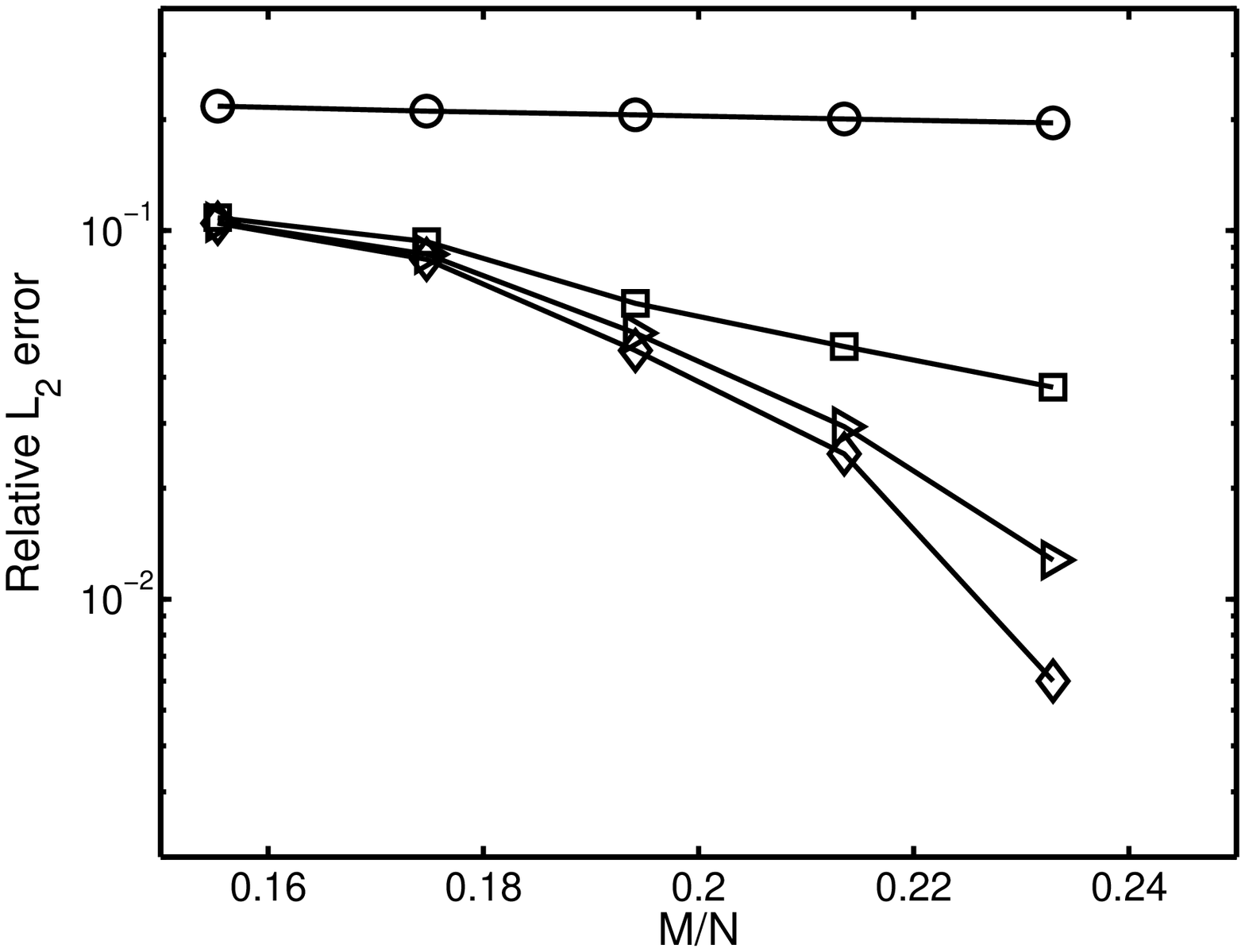}
\caption{Results for the example high-dimensional function. (Left) Comparison 
with $\ell_1$ methods.  ``$\circ$": standard $\ell_1$, ``$\Box$": rotated 
$\ell_1$ with $1$ iterations, ``$\triangleright$": rotated $\ell_1$ with $2$ 
iterations; ``$\diamond$": rotated $\ell_1$ with $3$ iterations. (Right)
Comparison with OMP methods.  ``$\Box$": rotated OMP with $1$ iterations,
``$\triangleright$": rotated OMP with $2$ iterations; ``$\diamond$": rotated 
OMP with $3$ iterations. These calculations were performed with $d=100$ and the
number of unknowns $N=5151$.}
\label{fig:ex5_l1} \label{fig:ex5_omp}
\end{figure}

\begin{figure}[t]
\centering
\includegraphics[width=3in]{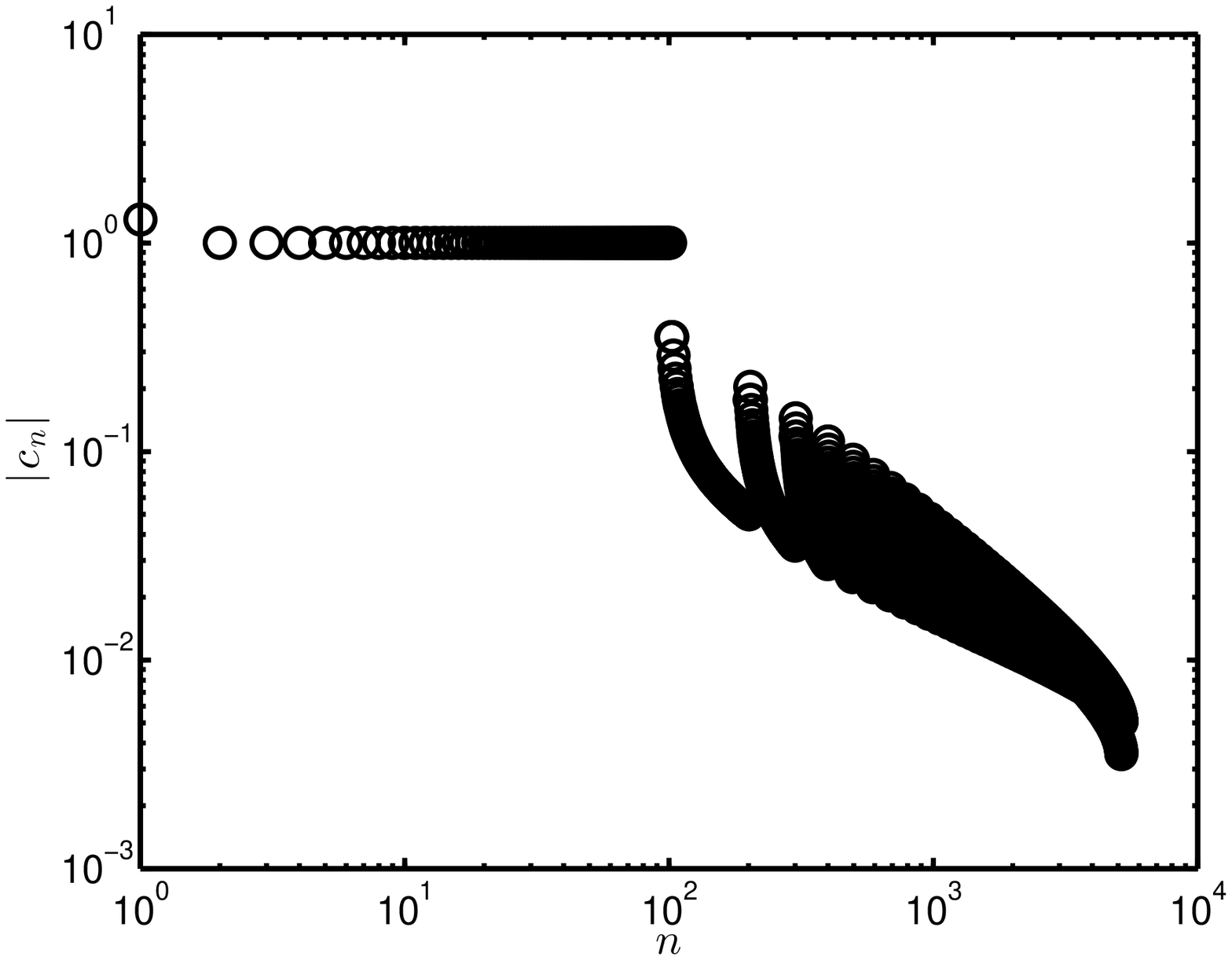}
\includegraphics[width=3in]{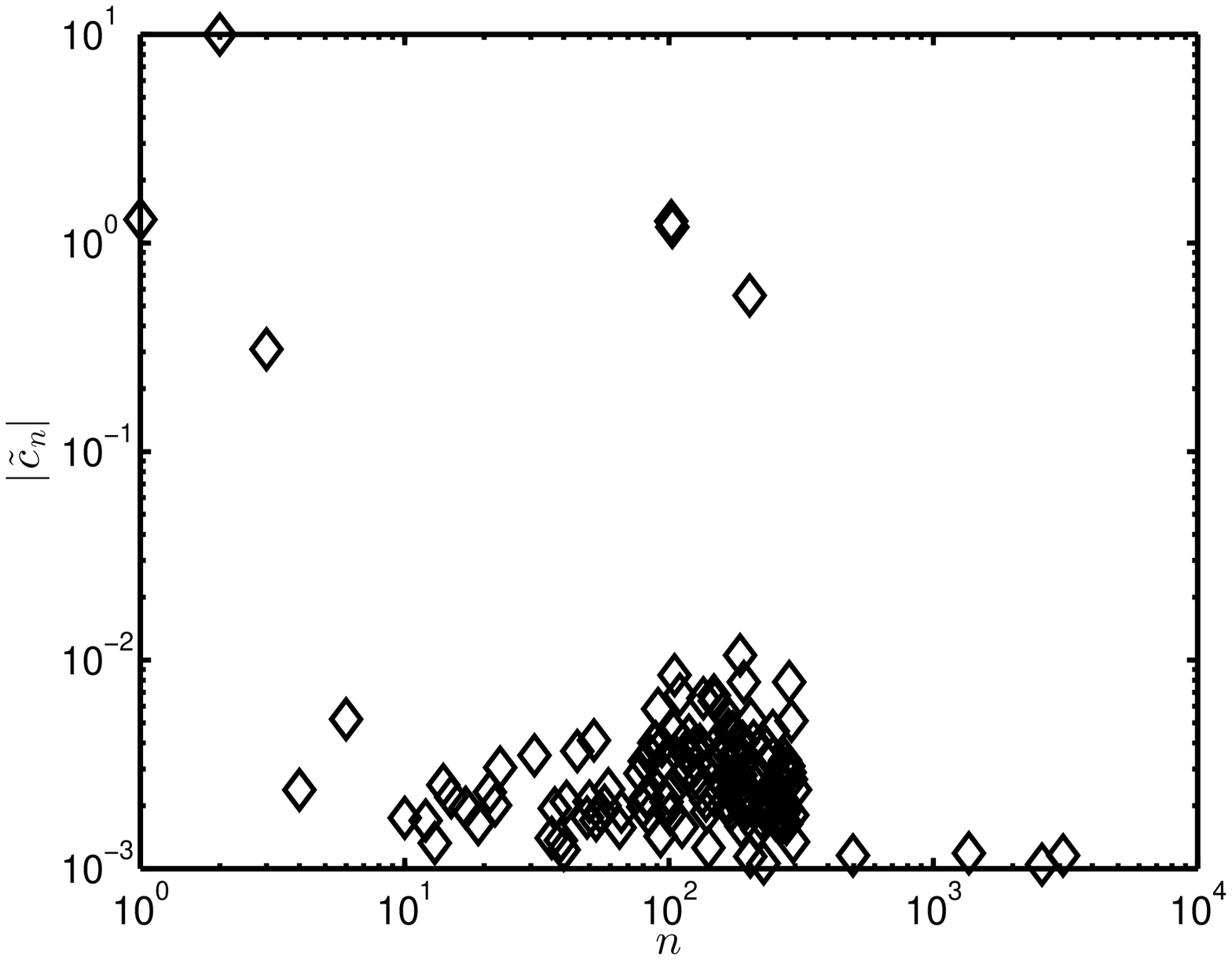}
\caption{Magnitude of the gPC coefficients for example high-dimensional 
function. (Left) Magnitude of $\bm c_n$. (Right) Magnitude of $\tilde{\bm c}_n$
of a randomly chosen replicate computed by $3$ iterated OMP method with $M=1200$
($M/N\approx 0.23$). These calculations were performed with dimension $d=100$ 
and the number of unknowns $N=5151$.}
\label{fig:ex5_coef}
\end{figure}

The previous examples demonstrate the capability of our new method to solve 
moderately high-dimensional problems. In the last example, we illustrate its 
potential for dealing with higher-dimensional problems. Specially, we select a
function similar to the first example (Sec.~\ref{subsec:ex1}) but with much 
higher dimensionality:
\begin{equation}
u(\bx) = \sum_{i=1}^d \xi_i + 0.25\left(\sum_{i=1}^d \xi_i/\sqrt{i}\right)^2,
  \quad d=100.
\end{equation}
The total number of basis functions for this example is $N=5151$. The relative 
error is computed with a level-$3$ sparse grid method, hence the numerical 
integrals are exact. The results are presented in Fig.~\ref{fig:ex5_l1}. As
before, our iterative rotation approach out-performs the existing $\ell_1$ and 
OMP methods. A comparison of $\bm c$ and $\tilde{\bm c}$ is presented in
Fig.~\ref{fig:ex5_coef} and it shows the enhancement of the sparsity by
the iterative rotation method.

For general high-dimensional problems, simply truncating the gPC expansion up to
a certain order is not efficient because the number of basis function will be 
very large. For example, in this test, $P=2$ requires $5151$ basis functions. 
Under such conditions even a small $M/N=0.2$ needs $1030$ samples, which can be
difficult in practical problems when the computational model is very costly.
Hence, a good approach for high-dimensional problems is to integrate our 
iterative rotation method with an adaptive method to reduce $N$; e.g., adaptive 
basis selection \cite{JakemanES14} or an ANOVA method \cite{YangCLK12}).

\subsection{Accuracy of computing the expansion coefficients $\bm c$}
\label{subsec:accuracy}
\begin{figure}[t]
\centering
\includegraphics[width=3in]{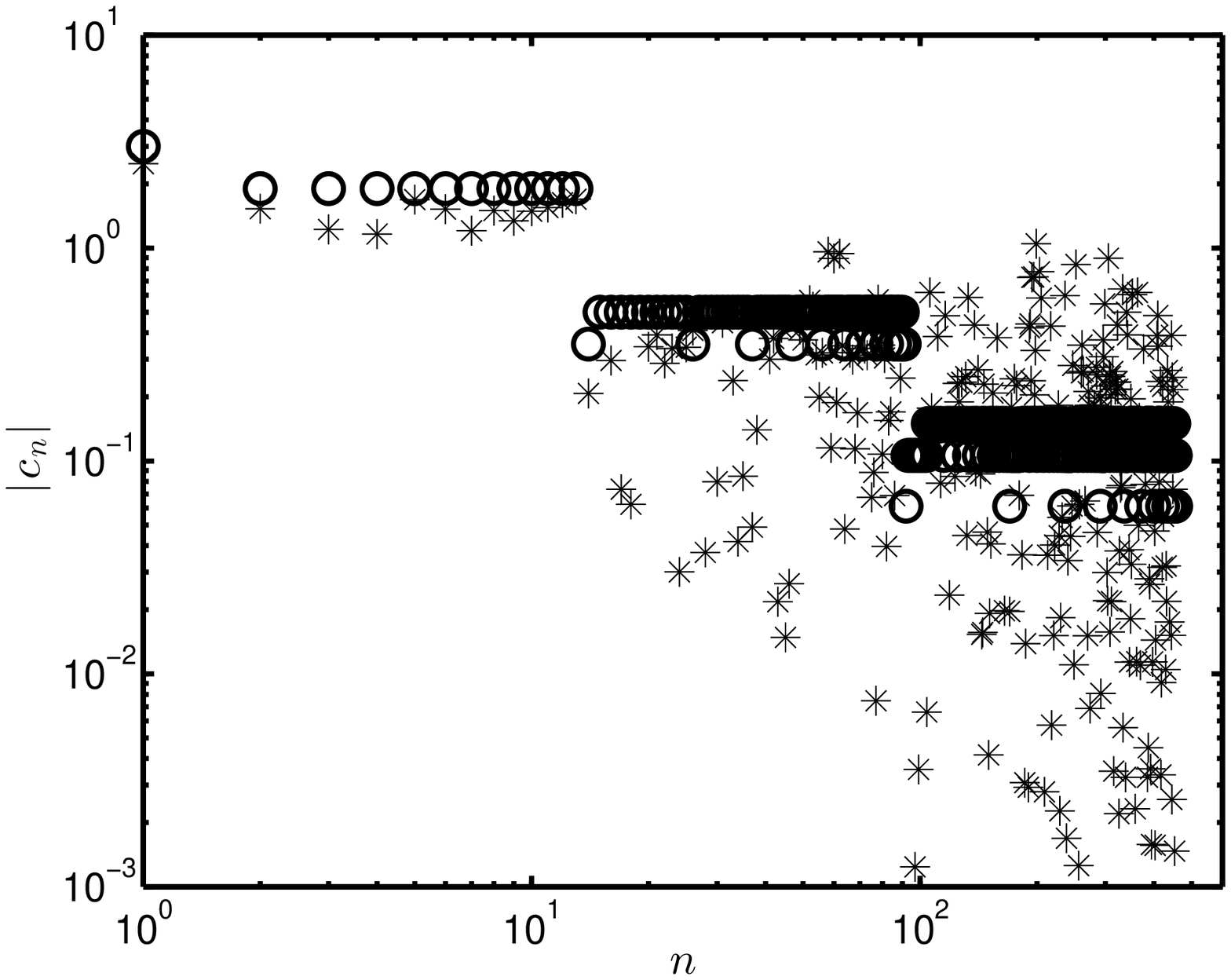}~
\includegraphics[width=3in]{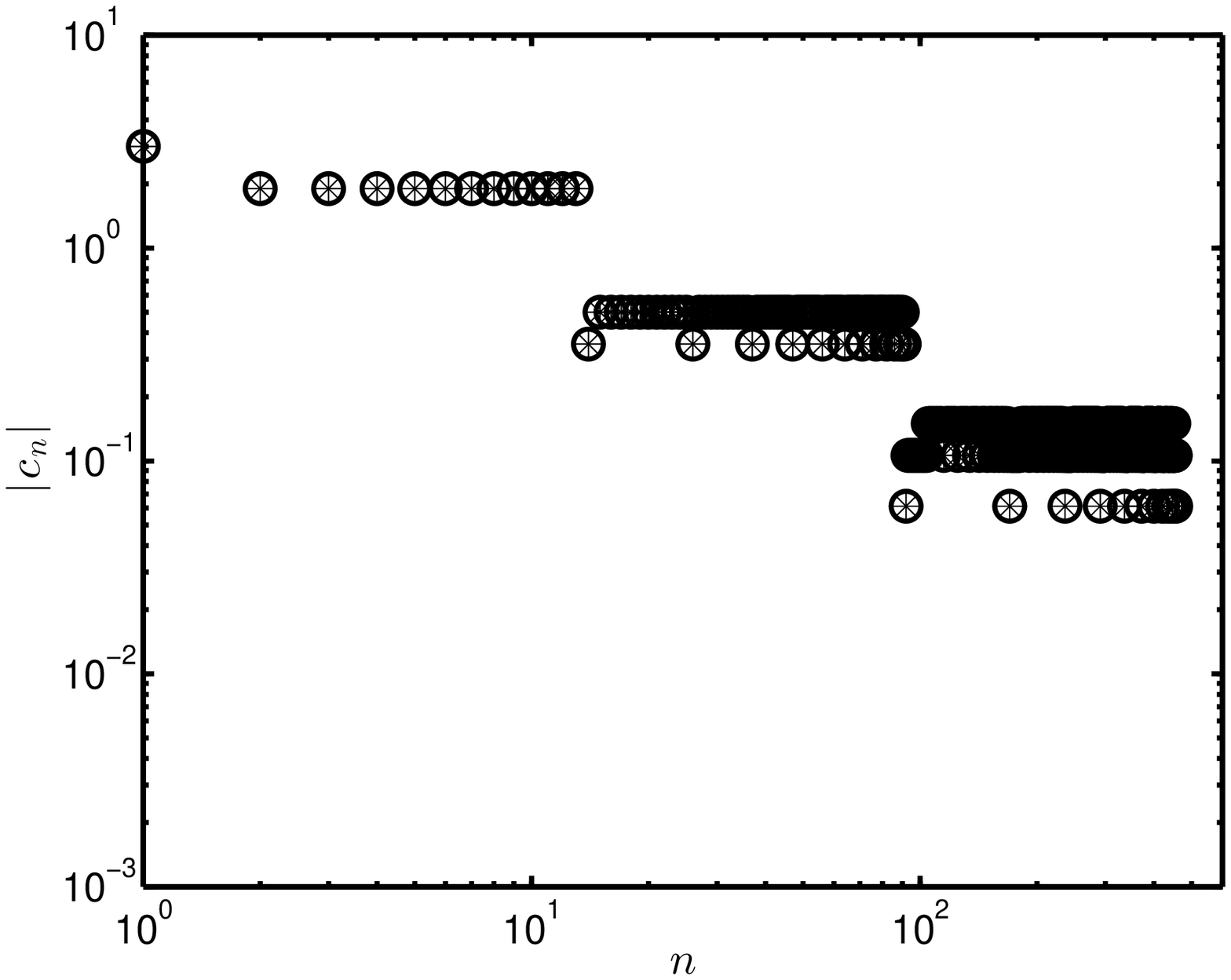}
\caption{Comparison of gPC coefficients for the example function with equally
important random  variables (Sec.~\ref{subsec:ex1}). (Left) Coefficients 
calculated by the standard $\ell_1$ method. ``$\circ$": exact $|c_i|$; ``$*$": 
$|c_i|$ by standard $\ell_1$ method. (Right) Coefficients calculated by our new 
iteratively rotated $\ell_1$ method with $L=9$ iterations.}
\label{fig:ex1_comp}
\end{figure}

\begin{figure}[t]
\centering
\includegraphics[width=3in]{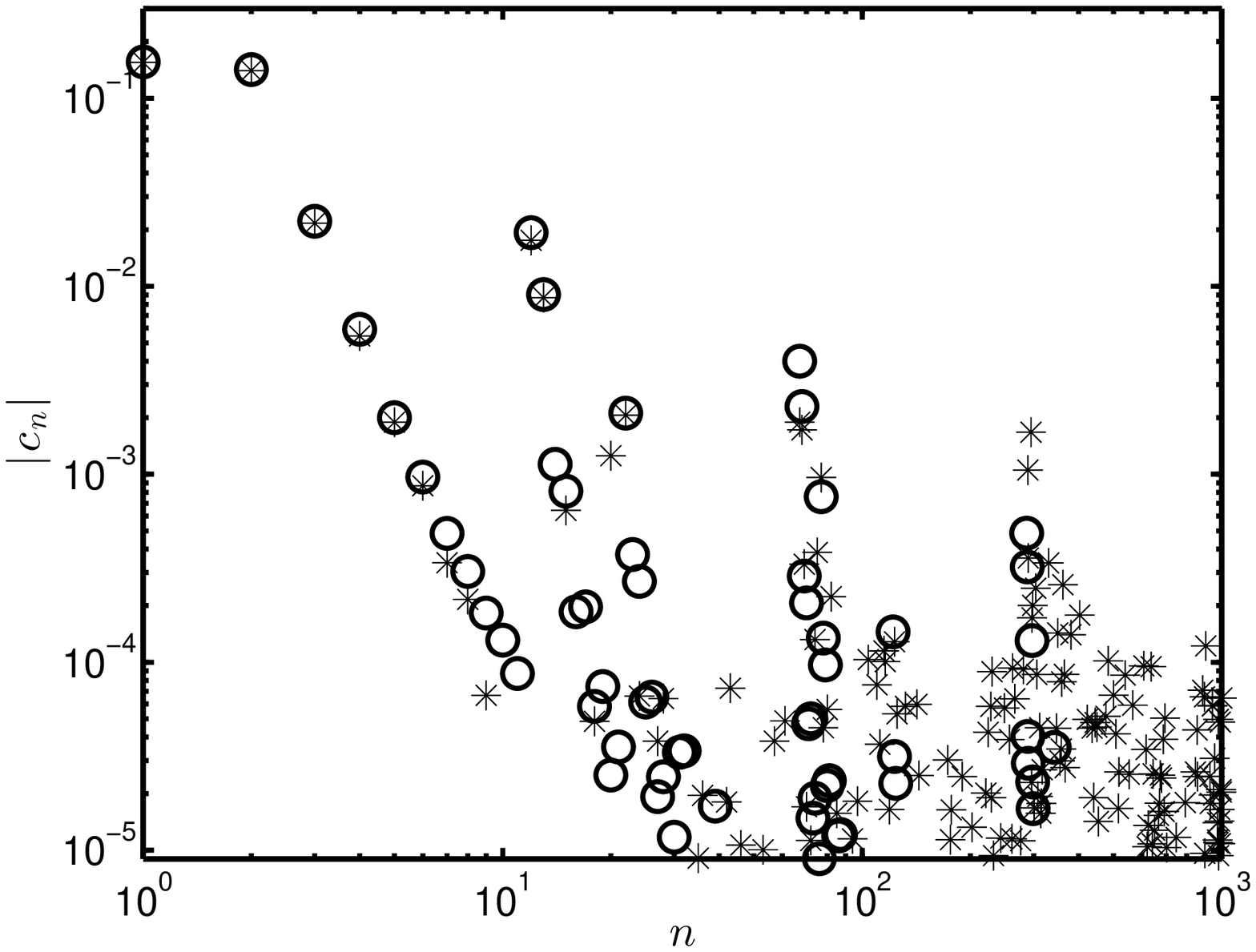}~
\includegraphics[width=3in]{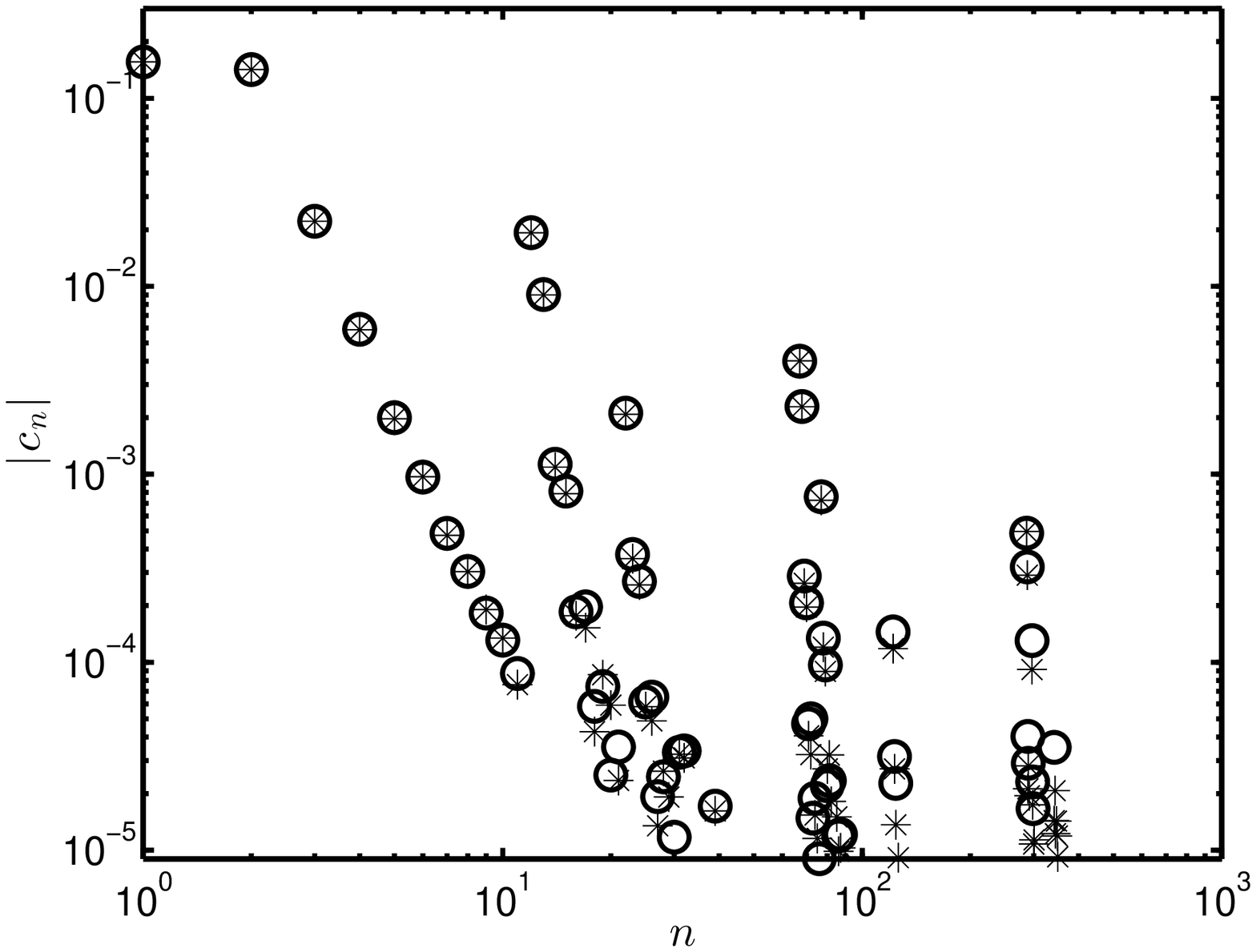}
\caption{Comparison of gPC coefficients for solutions of the KdV equation
(Sec.~\ref{subsec:ex4}). (Left) Coefficients calculated with the standard 
$\ell_1$ method. ``$\circ$": ``exact" $|c_i|$; ``$*$": $|c_i|$ by standard 
$\ell_1$ method. (Right) Coefficients calculated by our new iteratively rotated
re-weighted $\ell_1$ method (right). Only $|c_i|>10^{-5}$ are presented.}
\label{fig:ex4_comp}
\end{figure}

In many applications, gPC expansions are also used to study the sensitivity of 
the quantity of interest to the input random variables $\bx$. In order to 
perform this analysis, we need to transform the 
$u_g(\bm\eta)=\sum_{n=1}^N\tilde c_n\psi_n(\bm\eta) $ back to the original 
variables $u_g(\bx)=\sum_{n=1}^N c_n\psi_n(\bx)$. This transformation can be 
accomplished through inversion of $\tensor A$ in $\bm\eta=\tensor A\bx$. In 
Figures \ref{fig:ex1_comp} and \ref{fig:ex4_comp}, we present the coefficients 
of $u_g(\bx)$ in examples \ref{subsec:ex1} and \ref{subsec:ex4}, respectively. 
In both figures we randomly choose one test from the $100$ replicates. In
Figure \ref{fig:ex1_comp}, we select a result with $M=180$ ($M/N\approx 0.4$).
Using the standard $\ell_1$ method (left) gives very inaccurate results for
$\bm c$. However, Figure \ref{fig:ex1_comp} (right) shows that the iterative
rotation method with $L=9$ iterations gives much more accurate results for 
$\bm c$. We observe the same behavior in Figure \ref{fig:ex4_comp}, where we
chose a test with $M=120$ ($M/N\approx 0.12$) for the example KdV equation 
(Sec.~\ref{subsec:ex4}). In order to make this figure legible, we only present
those $c_i$ with absolute value larger than $10^{-5}$. This example demonstrates
that coefficients $c_n$ with magnitude larger than $10^{-3}$ are computed 
accurately by the combined iterative rotation and re-weighted $\ell_1$ method 
while the standard $\ell_1$ obtained significantly less accurate $c_n$. This 
difference is more distinct for $|c_n|\in [10^{-4},10^{-3}]$; in this range our
new method compute $c_n$ much more accurate than the standard $\ell_1$ method. 
In the lower right corner of the left plot, the standard $\ell_1$ method yields
many $c_n$ which should not appear in that area. As a comparison, we do not see
such $c_n$ calculated with the new method.

%% file: conclusion.tex
\section{Conclusions}
\label{sec:conclusion}

In this paper, we { extend our previous work \cite{LeiYZLB14}} 
and have introduced a compressive sensing-based gPC method to
increase the sparsity and accuracy of Hermite polynomial expansion with 
iterative rotations. Similar to the active subspace method 
\cite{ConstantineDW14,LeiYZLB14}, the rotation is decided by seeking the directions of 
maximum variation for the quantity of interest. Our current numerical examples 
are intended to demonstrate the ability of the method to increase the sparsity 
and accuracy of the gPC expansion; therefore, the quantity of interest only 
relies on the random variables. It is also possible to include the physical 
variables in the basis functions, i.e., 
$u(\bm x;\bx)=\sum_n c_n \psi_n(\bm x;\bx)$ (e.g, \cite{KaragiannisBL15}), our
future work will explore how this new method may help to increase the sparsity
in such cases.

We have demonstrated the method for $\ell_1$ minimization and OMP methods but it
can also be integrated with other compressive sensing methods. In particular, 
future work will investigate the integration of our new methods with advanced 
sampling strategies (e.g., \cite{HamptonD15}), adaptive basis selection method 
(e.g., \cite{JakemanES14}), Bayesian compressive sensing method 
(e.g.,\cite{KaragiannisBL15}), etc. These advanced strategies are particularly 
important for high-dimensional problems.

With this method, we will also be able to construct an accurate surrogate model
of the quantity of interest with limited data. Surrogate models are specifically
useful for the problems where the experiments or simulations are very costly.
This surrogate model can be used to study the sensitivity of the parameters and
is very useful in inverse problems based on Bayesian framework. Our new method
requires fewer output data to construct such surrogate models, which can be a 
great savings of experimental or computational resources.

Finally, we highlight three additional areas of future work for improving the
new method. First, it is currently only suitable for Hermite polynomial 
expansions. Second, the new method requires a formal numerical analysis to 
assess convergence behavior and determine specific terminating criteria. 
Finally, there are likely more optimal iterative rotation strategies that can be
applied to Hermite polynomial or other expansions. One possible direction of 
work is to design a suitable objective function and consider this problem from
an optimization point of view. All of these questions will be addressed in our
future work.

\section*{Acknowledgments}
This work is supported by the U.S.~Department of Energy, Office of Science,
Office of Advanced Scientific-Computing Research as part of the Collaboratory
on Mathematics for Mesoscopic Modeling of Materials (CM4). Pacific Northwest
National Laboratory is operated by Battelle for the DOE under Contract 
DE-AC05-76RL01830.

%% file: append.tex
\section*{Appendix}
\renewcommand{\theequation}{A-\arabic{equation}}

Here we provide the details of computing the elements $(K_{ij})_{kl}$ in Eq.~\eqref{eq:kernel}.
Notice that for univariate normalized Hermite polynomials, 
\begin{equation}
\psi_n(\xi)'=\sqrt{n}\psi_{n-1}(\xi),\quad n\in\mathbb{N}\cup\{0\},
\end{equation}
where we set $\psi_{-1}(\xi)=0$ for simplicity. Therefore, we have
\begin{equation}
\mexp{\psi_i(\xi)'\psi_j(\xi)}=
\int_{\mathbb{R}}\psi_i(\xi)'\psi_j(\xi)\rho(\xi)\dif\xi=\sqrt{i}\delta_{i-1 j},
\end{equation}
where $\rho(\xi)$ is the PDF of a standard Gaussian random variable. For a
multi-index $\ba=(\alpha_1,\alpha_2,\cdots,\alpha_d), \alpha_i\in
\mathbb{N}\cup\{0\}$, and basis function 
$\psi_{\ba}(\bx)=\psi_{\alpha_1}(\xi_1)\psi_{\alpha_2}(\xi_2)\cdots
\psi_{\alpha_d}(\xi_d)$,
\begin{equation}
\dfrac{\partial}{\partial\xi_i}\psi_{\ba}(\bx)=\psi_{\alpha_i}(\xi_i)'
\prod_{\substack{m=1\\ m\neq i}}^d \psi_{\alpha_m}(\xi_m).
\end{equation}
Hence, given two different multi-indices
$\ba_k=((\alpha_k)_{_1},(\alpha_k)_{_2},\cdots,(\alpha_k)_{_d})$ and 
$\ba_l=((\alpha_l)_{_1},(\alpha_l)_{_2},\cdots,(\alpha_l)_{_d})$, the
corresponding entry of matrix $\tensor K_{ij}$ is
\begin{equation}
\begin{aligned}
(K_{ij})_{kl} & =  \mexp{\dfrac{\partial\psi_{\ba_k}(\bx)}{\partial\xi_i}\cdot
\dfrac{\partial\psi_{\ba_l}(\bx)}{\partial\xi_j}} \\
& =  \mexp{
\left(\psi_{(\alpha_k)_{_i}}(\xi_i)'\prod_{\substack{m=1\\ m\neq i}}^d
    \psi_{(\alpha_k)_{_m}}(\xi_m)\right)\cdot
\left(\psi_{(\alpha_l)_{_j}}(\xi_j)'\prod_{\substack{m=1\\ m\neq j}}^d
    \psi_{(\alpha_l)_{_m}}(\xi_{m})\right)} \\
& = \sqrt{(\alpha_k)_{_i}(\alpha_l)_{_j}}
\delta_{(\alpha_k)_{_i}-1(\alpha_l)_{_i}}
\delta_{(\alpha_k)_{_j}(\alpha_l)_{_j}-1}\cdot
\prod_{\substack{m=1\\ m\neq i,m\neq j}}\delta_{(\alpha_k)_{_m}(\alpha_l)_{_m}}.
\end{aligned}
\end{equation}